\documentclass[11pt]{article}

\usepackage{amssymb}

\font\tenbb=msbm10
\def\cC{\hbox{\tenbb C}}
\def\rR{\hbox{\tenbb R}}
\def\nN{\hbox{\tenbb N}}
\def\qQ{\hbox{\tenbb Q}}

\def\un{\underline}
\def\td{\bigtriangledown}

\def\di{\displaystyle}

\def\re{\mbox{\rm Re}}
\def\im{\mbox{\rm Im}}

\newtheorem{thm}{Theorem}
\newtheorem{defi}{Definition}
\newtheorem{lem}{Lemma}

\newtheorem{cor}{Corollary}
\newtheorem{rema}{Remark}
\newtheorem{proper}{Property}
\newtheorem{sta}{Statement}
\newtheorem{ass}{Assumption}

\begin{document}
\baselineskip 6mm
\def\td{\bigtriangledown}
\titlepage

\title{Scale calculus and the Schr\"odinger equation}
\author{Jacky CRESSON}

\maketitle

\begin{center}
Equipe de Math\'ematiques de Besan\c{c}on,

CNRS-UMR 6623, Equipe Alg\`ebre et Th\'eorie des nombres,

16 route de Gray,

25030 Besan\c{c}on cedex, France

e-mail: cresson@math.univ-fcomte.fr
\end{center}

\begin{abstract}
This paper is twofold. In a first part, we extend the classical differential calculus to
continuous non differentiable functions by developping the notion of scale calculus.
The scale calculus is based on a new approach of continuous non differentiable functions
by constructing a one parameter family of differentiable functions $f(t,\epsilon )$ such that $f(t,\epsilon )
\rightarrow f(t)$ when $\epsilon$ goes to zero. This lead to several new notions
as representation, fractal functions and $\epsilon$-differentiability. The basic
objets of the scale calculus are left and right quantum
operators and the scale operator which generalize the classical derivative.
We then discuss some algebraic properties of these operators.
We define a natural bialgebra, called quantum bialgebra,
associated to them. Finally, we discuss a convenient geometric object associated
to our study. In a second part, we define a first quantization procedure of classical
mechanics following the scale relativity theory developped by Nottale. We obtain a
non linear Schr\"odinger equation via the classical Newton's equation
of dynamics using the scale operator. Under special assumptions we recover the
classical Schr\"odinger equation and we discuss the relevance of these assumptions.
\end{abstract}


\newpage
\tableofcontents

\newpage

\section*{Introduction}

The origin of the fundamental incompatibility between {\it quantum mechanics} and
Einstein's {\it general relativity} lies in the microscopic geometric structure of
space-time. As pointed out by Greene \cite{gr}, Feynman \cite{fey}, Cohen-Tannoudji and Spiro (\cite{cospi},
p.131) and others, space-time is no more a {\it differentiable} manifold at the {\it atomic}
scale, contrary to the assumption of general relativity.\\

>From this fact, at least two theory have been constructed : \\

- the {\it string theory}, which implies a dimensional extension of space-time by
allowing closed dimension at the Planck scale;\\

- the {\it scale relativity theory} developped by Nottale \cite{no1}, which gives up the
Einstein's assumption of the differentiability of space-time by considering what he calls a
{\it fractal space-time}, which can be interpreted as a scale dependant non differentiable
manifold. He then extends einstein's relativity principle to scale, and developp the
{\it scale relativity principle}.\\

In this article, we explore this second alternative.\\

Nottale \cite{no1} has studied what are the consequences of the abandon of the differentiability
of space-time. This problem is difficult, in particular because the mathematical foundations of
such a theory are not yet constructed. For example, Nottale asserts that there exists an
``infinity of geodesics on a fractal space-time". This sentence is difficult to understand because
we don't know what is a fractal space-time\footnote{We refer to \cite{c4} for a first
definition of a fractal manifold and a discussion of the special scale relativity theory.} and even if we identify this set to a non differentiable
manifold, we don't know what is the sense of ``geodesic". As a consequence, we restrict
our attention to a far simpler problem, namely the consequences of the loss of
differentiability of a given trajectory. A first approach is to consider that only
{\it space} is a non differentiable manifold and to take the variable $t$ as an absolute
variable. As a consequence, trajectories of quantum particles are non differentiable
curves parametrized by the time variable $t$. In this case, we have at least the following
two consequences :\\

i) By Lebesgues theorem \cite{tri}, the {\it lenght} of
a non differentiable curve $\Gamma$ is {\it infinite}. What does it means from the physical view-point ?
That given a parameter $\epsilon >0$, which has the sense of a resolution, the lenght $L_{\epsilon}$
of the curve constructed by connecting small arcs of lenght $\epsilon$ on points of $\Gamma$,
goes to infinity when $\epsilon$ goes to zero. As a consequence, the role of $\epsilon$ is now
fundamental as $L_{\epsilon}$ loes any sense when $\epsilon$ going to zero,
contrary to usual differentiable curves where $\epsilon$ is only a parameter
of precision and $L_{\epsilon}$ goes to a fixed constant $L$.

Nottale then introduce {\it fractal functions}, which are
resolution (or scale) dependant functions $f(x,\epsilon )$, which converge to non differentiable
functions, and a ({\it renormalization like}) differential equation satisfied by $L_{\epsilon}$,
called a {\it scale law}, which gives the behaviour of $f(x,\epsilon )$ when $\epsilon$ goes
to zero.\\

ii) The derivative along the curve has no sense. Nottale introduce a complex operator,
that he calls the {\it scale derivative}. It takes into account the mean-backward and
mean-forward derivative along the curve.\\

Using these tools, he gives an informal derivation of the {\it Schr\"odinger equation} from
the classical {\it Newtonian equation of dynamics}, via a quantization
procedure which follows from an extension of Einstein's relativity principle called
the {\it scale relativity principle}.\\

In this paper, we developp a mathematical framework in which we can explicit the
quantization procedure, which we call the {\it scale quantization procedure}.\\

The plan of the paper is as follows :\\

In part 1, we define a natural extension
of Leibniz differential calculus which can be used on non differentiable functions
in order to precise points i) and ii).
We introduce the {\it Scale calculus}, which formalize the concept of
$\epsilon$-differentiability. In particular, we define an operator called
{\it scale difference operator}, which is the rigourous mathematical counterpart of
Nottale's scale derivative.\\

In part 2, we define the scale
quantization procedure. We give a precise definition of the quantization map,
which allows us to associate to the classical Newtonian equation of dynamics a
quantized analog. This analog has the form of a generalized non linear Schr\"odinger
equation. We then discuss how to obtain the classical Schr\"odinger equation.

\newpage
\part{Scale calculus}

\section{Introduction}

Non differentiable functions, and more generally non differentiable manifolds,
become more and more important in many part of mathematics and physics, like
brownian motion \cite{ei}, quantum mechanical path by Feynman and Hibbs \cite{fh}. Despite
many works, our understanding of non differentiable functions is insatisfactory.

A great deal of efforts have been devoted to generalize, as long as possible, the
classical differential calculus of Leibniz and Newton. This lead to different kind
of fractional calculus (Riemann, Liouville, Weyl,...). All these fractional
calculus are based on a pure analytic generalization of the
Cauchy formula. As a consequence, and despite their intrinsic interest, they
are difficult to interpret (in particular, from the geometrical view
point).\\

The aim of this paper is to introduce a set of ideas, coming from physics, in order to
renew our approach to non differentiable functions.\\

In physic, the non differentiability is not studied for itself. In the contrary,
this is the effect of non differentiability with respect to a differentiable
model wich is looking for. For example, as explained by Greene (\cite{gr}, chap.5), in {\it superstring} theory
one is lead to a new vision of spacetime because the spacetime, at the scale of
atoms can not be considered differentiable at all, as in the general relativity scheme.
Here, one focus on fluctuations with respect to a differentiable character. \\

Moreover, one usually doesn't have access via measurement, to the non differentiable object
(function, manifold), but to a almost everywhere differentiable model of it. This explain
also the previous remark : by measure we obtain a differentiable model, and we must see,
when the precision of measure increases, if there is no two strong fluctuations with
respect to this model. For a differentiable process, the fluctuations decrease.
In the case of a non differentiable process, one expect larger and larger
fluctuations.\\

This point of view lead to several new concepts : representation of continuous non differentiable
functions, fractals functions and $\epsilon$-differentiability. The idea is to associate to each
continuous non differentiable functions $f(t)$, a one parameter family of differentiable function
$f(t,\epsilon )$ such that $f(t,\epsilon ) \rightarrow f(t)$ when $\epsilon$ goes to zero. The next
step is to give a usefull criterion which say that $f(t,\epsilon )$ is a ``good" model for
the function $f$. An important notion is then the minimal resolution wich is, more or less,
the precision under which, one can't use a
differentiable model without laking many features of the underlying function $f$. We then define quantum
derivatives and the scale derivative, which reflect the non differentiable character of
the underlying function $f$. \\

We discuss algebraic properties of quantum derivatives, and
the scale derivative. These operators act on the set of continuous real valued functions,
denoted by $C^0$. They are first introduced in \cite{bc1} in order to discuss the
derivation of the Schr\"odinger's equation from the classical Newton's equation of
dynamics using Nottale's scale relativity theory \cite{no1}. \\

In this paper, we construct a natural structure of bialgebra, called quantum
bialgebra, using specific properties of quantum derivatives. The quantum
algebra can be considered as a small deformation of a classical Hopf
algebra. Although we are close to problem related to quantum groups and
quasitriangular Hopf algebras introduced by Drinfeld \cite{dr1}, we stress that
quantum bialgebra are new.

There is no natural bialgebra structure associated to scale derivative. This follows
from the lack of a natural composition rule for scale operators that we define in
this paper. \\

We also discuss the natural geometric object associated to a continuous non differentiable
function. It turns out that this geometric object is, in the simplest case, the product of
a differentiable curve $\Gamma$ by a two points set $A=\{ a,b\}$, so $M\times A$.
This is a simple example of a noncommutative space studied by Connes \cite{co},
as a preliminary to his noncommutative model of the standard model.

\section{About non differentiable functions}

The aim of this paragraph is to develop a rigourous mathematical understanding of the idea
of fluctuation with respect to a differentiable model for continuous non differentiable functions.

\subsection{Representation theory and fractals functions}

Here, we introduce two dual notions : representation of non differentiable functions and
fractals functions. Representation are well suited to develop a mathematical understanding
of non differentiable function based on differentiable functions. Fractals functions takes
their origin in physical problems. \\

In the following, we denote by $C^0$ the set of continuous real valued functions, and by $C^1$
the set of differentiable real valued functions.

\begin{defi}
\label{d1}
Let $f\in C^0$, a representation of $f$ is a one parameter family of operator $S_{\epsilon}$,
defined by
\begin{equation}
\label{repre}
S_{\epsilon} :
\left .
\begin{array}{lll}
C^0 & \longrightarrow & C^1 \\
f & \longmapsto & S_{\epsilon} (f)=f_{\epsilon} ,
\end{array}
\right .
\end{equation}
and such that the differentiable functions $f_{\epsilon} \in C^1$ converge, in
$C^0$-topology, toward $f$ when $\epsilon$ goes to zero.
\end{defi}

A basic example is obtain by approximating $f$ by mean functions $f_{\epsilon}$ defined by
$f_{\epsilon} (t) =(1/2\epsilon )\di\int_{t-\epsilon}^{t+\epsilon} f(s) ds$. \\

More generally, we can consider a {\it smoothing function} $\Phi (s;t,\epsilon )$ depending
on two parameters, $t$ and $\epsilon$, satisfying the normalization
\begin{equation}
\di\int_{-\infty}^{\infty} \Phi (s;t,\epsilon ) ds =1 .
\end{equation}
For any continuous function, we define a representation by
\begin{equation}
f_{\epsilon} (t) =\di\int_{-\infty}^{\infty} \Phi(s;t,\epsilon ) f(s) ds .
\end{equation}

In practice, we have never access to $f$. An idea is to define $f$ via a family of functions
which are not functionally dependent of $f$ like in (\ref{repre}). We are lead to
the notion of {\it fractal functions}, first introduced by Nottale \cite{no1}
(see also \cite{c2}).

\begin{defi}
\label{d2}
A fractal function is a parametrized function of the form $f(t,\epsilon )$, depending on
$\epsilon >0$, such that : \\

i) For all $\epsilon >0$, the function $f(t,\epsilon )\in C^1$, except at a finite number of points,\\

ii) There exists an everywhere non differentiable function $f(t)$ such that $f(t,\epsilon )$
converge to $f(t)$ when $\epsilon$ goes to zero.
\end{defi}

The main difference between definition \ref{d1} and \ref{d2} is that for fractal functions,
one usually doesn't know an explicit form of the limit function $f$,
We only require an existence result.
Moreover, the set of functions $f(t,\epsilon )$ doesn't refer to the limit $f$ in its
definition, which is closest to the  {\it measurement} process in physical experiment. \\

\subsubsection{Examples of fractal functions : Nottale's functions and iteration
of affine systems}
\label{nota}

The basic example of fractal functions is {\it Nottale's functions} introduced in
\cite{no1} :\\

For all $\epsilon >0$, and for all $0<\mu <\epsilon$,
\begin{equation}
\label{nottale}
x(t,\epsilon )= \di\int \Phi_{\epsilon ,\mu } (t,y) x(y,\mu ) dy ,
\end{equation}
where $\Phi_{\epsilon ,\mu} (x,y)$ is a differentiable function such that
\begin{equation}
\label{smooth}
\int_{-\infty}^{\infty} \Phi_{\epsilon ,\mu} (x,y) dy =1 , \
\forall \, x\in {\bf R} ,
\end{equation}
called a smoothing function.

\begin{defi}
Let $\Phi_{\epsilon ,\mu}$ be a smoothing function satisfying (\ref{smooth}).
We denote by ${\cal N} (\Phi_{\epsilon ,\mu } )$, and we call Nottale's set associated to
$\Phi_{\epsilon ,\mu}$, the set of functions defined by (\ref{nottale}).
\end{defi}

We refer to \cite{c2} for basic properties of this set of functions, in particular for a
usefull equivalence relation.\\

An interesting example of fractal functions for which the limit function is not explicit
is given by iteration of {\it affine systems} \cite{tri}. \\

An affine map in $\rR^2$, with a coordinates system $(x,y)$, is a map of the form
\begin{equation}
\label{affine}
F
\left (
\begin{array}{c}
x \\
y
\end{array}
\right )
= M
\left (
\begin{array}{c}
x \\
y
\end{array}
\right )
+ T ,
\end{equation}
where $M$ is a $2\times 2$ matrix and $T$ is a translation vector.\\

An {\it affine system} is given by : \\

i) a positive integer $N\geq 2$,

ii) $N+1$ points of $\rR^2$, $A=A_1$, $\dots$, $A_{N+1} =B$. We denote
$A_i =(x_i ,y_i)$, $x_1 =a$, $x_{N+1} =b$. We assume that
$$a=x_1 < x_2 <\dots < x_{N+1} =b .$$
iii) $N$ affine map $F_1 ,\dots ,F_N$ such that
$$F_i (AB)=A_i A_{i+1} .$$
We denote by $F$ the map defined by
\begin{equation}
F(E) =\di\bigcup_{i=1}^N F_i (E) .
\end{equation}
Let $z_0$ be the affine function on $[a,b]$ whose graph $\Gamma_0$ is the segment $AB$.
The image $\Gamma_1 =F(\Gamma )$ is the graph of a continuous function, which is affine.
For all $n$, we define the continuous function $z_n$ whose graph is $\Gamma_n =
F(\Gamma_{n-1} )$. Following (\cite{tri}, p.175), the sequence of functions $(z_n )_{n\in \nN}$ converge uniformly to
a continuous function $z_{\infty}$ such that $F(\Gamma_{\infty} )=\Gamma_{\infty}$.

\subsection{$\epsilon$-differentiability and minimal resolution}

A basic properties of differentiable functions is that the quantities
\begin{equation}
\td_+^{\epsilon} f (t) = \di {f(t+\epsilon )-f(t) \over \epsilon} ,\ \mbox{and}\ \td_-^{\epsilon} f(t) =
\di {f(t)-f(t-\epsilon )\over \epsilon} ,
\end{equation}
keep sense when $\epsilon$ goes to zero and are equal. \\

As a consequence, the following quantity
\begin{equation}
a_{\epsilon} f(t) =\left | \di {f(t+\epsilon )+f(t-\epsilon )-2 f(t) \over \epsilon} \right | ,
\end{equation}
converges to zero when $\epsilon$ goes to zero. \\

The underlying idea is that the two representations of a function $f$, given by the
forward and backward mean function, defined as
\begin{equation}
f_{\epsilon}^+ (t)=(1/2\epsilon) \di\int_t^{t+\epsilon} f(s) ds
\ \ \ \mbox{\rm and}\ \ \
f_{\epsilon}^- (t)=(1/2\epsilon ) \di\int_{t-\epsilon}^t f(s) ds,
\end{equation}
respectively, must have derivatives which coincide when $\epsilon$ goes to zero. \\

This remark allows us to introduce the following notion of $\epsilon$-$h$-differentiability :

\begin{defi}
Let $h>0$ be a given real number. A function $f\in C^0$ is says to be
$\epsilon$-$h$-differentiable at point $t$, if
\begin{equation}
a_{\epsilon} f(t) <h .
\end{equation}
\end{defi}

We can detect the non differentiable character of a function by investigating its
$\epsilon$-$h$ differentiability. Precisely, we define the notion of minimal resolution :

\begin{defi}
Let $h>0$ be a given real number and $f\in C^0$. The $h$-minimal resolution of $f$ at point $t$,
denoted $\epsilon(f,h)(t)$ is defined as $\inf_{\epsilon} \{ a_{\epsilon} f(t) <h\}$.
\end{defi}

Of course, if for a given $h$, the $h$-minimal resolution is non zero, then $f$ is
non differentiable. \\

For all $\alpha \in ]0,1[$, we denote by $C^{\alpha}$ the set of continuous real
valued functions, defined on $[0,1]$ such that the quantity
\begin{equation}
\mid f\mid_{\alpha} =\sup_{0\leq x\not= y\leq 1} \di {\mid f(x)-f(y)\mid \over \mid x-y\mid^{\alpha}} ,
\end{equation}
is finite (H\"olderian functions of order $\alpha$). Then, we have :

\begin{lem}
Let $0<\alpha <1$ and $f\in C^{\alpha}$. For all $t \in ]0,1[$, and all $h>0$, the $h$-minimal
resolution of $f$ at point $t$ satisfies
\begin{equation}
\epsilon (f,h)(t) \leq \di \left ( {h\over 2\mid f\mid_{\alpha} } \right ) ^{1/\alpha -1} .
\end{equation}
\end{lem}

{\bf Remark}. In this example, the minimal resolution depends on $\mid f\mid^{\alpha}$. As a
consequence, for a quantum mechanical path, we expect that the minimal resolution depends on
the momentum of the particle. This is indeed the case in Nottale's theory \cite{no1}, where the
minimal resolution is related to the De Broglie lenght of the particle. \\

A global order of $h$-minimal resolution can be defined.

\begin{defi}
Let $h>0$ be given and $f\in C^0$. The $h$-minimal resolution of $f$, denoted $\epsilon (f,h)$, is defined
by $\epsilon (f,h) =\sup_{t\in {\cal D}f} \epsilon (f,h)(t)$, where ${\cal D} f$ is the definition domain
of $f$.
\end{defi}

In this definition it is important to take the sup of $h$-minimal resolution of $f$ at point $t$. As a consequence,
if $f$ is differentiable on a small set of point, the $h$-minimal resolution is however non zero. \\

\begin{rema}
We have $\epsilon (f+c,h)=\epsilon (f,h)$ for all $c\in \rR$. But, we have
$\epsilon (\lambda f,h) \not= \epsilon (f,h)$ for all $\lambda \not=1$, contrary to the
case of \cite{bc1}. This inequality is related to the fact that a changing momentum induce
a change of regularity for the curve. A physical consequence, is that the minimal
resolution must depends on the momentum.
\end{rema}

The connexion to representation theory is done through the backward and forward mean functions,
$f_{\epsilon}^{\sigma}$, $\sigma=\pm$, introduced before. Indeed, backward and forward mean
functions are differentiable functions. They can be used as a classical representation of a given
continuous function $f$. However, if $f$ admits a non zero $h$ minimal resolution, this means that
these functions are not sufficient to capture the complete local behaviour of $f$ as long as
$\epsilon <\epsilon (f,h)$. Of course, this notion depends on $h$. In physical problems, the constant
$h$ must corresponds to a universal constraint, like Heisenberg constraint in quantum mechanics.

\subsection{Scale law}

For an everywhere non-differentiable function $f$, the lenght $L_{\epsilon}$ of the graph of the mean
function $f_{\epsilon}$ goes to infinity when $\epsilon$ goes to $0$\footnote{
Of course, this property can't be used as a definition of an everywhere non differentiable
function. We can find curve which are rectifiable with infinite lenght (see \cite{bc3}).}. We want to
quantify this divergence. A first idea is to find a differential equation which gives the
behaviour of $L_{\epsilon}$ with respect to $\epsilon$ like in (\cite{c2},\cite{bc3}).
However, this is difficult because $L_{\epsilon}$ is not differentiable with respect to $\epsilon$
in general. In this section, we define a less rigid definition of scale law and we discuss
its properties.

\subsubsection{Definition}

Let $f$ be an everywhere non differentiable continuous function on the interval $I=[0,1]$.
For $\epsilon >0$, we denote by $f_{\epsilon}$ the mean function associated to $f$ on $I$,
and ${\cal L}_{\epsilon}$ the lenght of its graph.

\begin{defi}
\label{scalaw}
We say that $f$ satisfies a scale law if there exist functions $l_{\epsilon} >0$ and
$L_{\epsilon} >0$ such that
\begin{equation}
l_{\epsilon} \leq {\cal L}_{\epsilon} \leq L_{\epsilon} ,
\end{equation}
satisfying
\begin{equation}
{\cal L}_{\epsilon} =O(l_{\epsilon} ) ,\ {\cal L}_{\epsilon} =O(L_{\epsilon} ) ,
\end{equation}
and for which there exist a function $E : \rR \rightarrow \rR$, such that
\begin{equation}
\di {d l_{\epsilon} \over d\ln \epsilon } =E(l_{\epsilon} ,\ln \epsilon ) , \
\di {d L_{\epsilon} \over d\ln \epsilon } =E(L_{\epsilon} ,\ln \epsilon ) .
\end{equation}
The function $E$ is called a scale law.
\end{defi}

Of course, when one knows a scale law of $f$, we deduce a speed of drift for
${\cal L}_{\epsilon}$.\\

Basic examples of scale laws are given by \\

i) For $a>0$, $E(x,t)=a$,\\

ii) For $b>0$, $E(x,t)=bx$.\\

For i), we obtain graph with logarithmic drift of order $\mid a \mid \ln (1/\epsilon )$.
For ii), we obtain classical power law drift of order $1/\epsilon^b$.\\

Using this notion, we are lead to two different kind of problems :\\

i) Let $E$ be a given function. Find the set of function ${\cal E} (E)$ such that for
all $f\in {\cal E}(E)$, a scale law of $f$ is $E$. \\

ii) Let ${\cal E}$ be a given set of functions. Find, if it exits, a scale law for
each $f\in {\cal E}$.\\

The first problem is equivalent to estimate a given class of functions by the speed of
drift of ${\cal L}_{\epsilon}$. This point is difficult and discussed in \cite{bc3}.\\

The second one is more natural and discussed in the following.

\subsubsection{Scale law of H\"olderian functions}

We first define an important class of continuous functions.

\begin{defi}
We denote by $H^{\alpha } (c,C)$ the set of real valued continuous functions $f$ such that
for all $\epsilon >0$ sufficiently small, and $\mid t-t' \mid <\epsilon$, we have
\begin{equation}
c\epsilon^{\alpha } \leq \mid f(t)-f(t') \mid \leq C\epsilon^{\alpha } .
\end{equation}
\end{defi}

The set $H^{\alpha} $ corresponds to continuous functions which are H\"oldet and
inverse H\"older of exponent $\alpha$.\\

{\bf Example}. Let $0<\alpha <1$ and $g(t)$ be the function of period $1$ defined on
$[0,1]$ by
\begin{equation}
g(t)=
\left \{
\begin{array}{lll}
2t , & \ \ \ & if\ 0\leq t\leq 1/2 , \\
2-2t, & & if\ 1/2 \leq t \leq 1 .
\end{array}
\right .
\end{equation}
The Knopp or Takagi function, defined by
\begin{equation}
K(t)=\di\sum_{n=0}^{\infty} 2^{-n\alpha } g(2^n t) ,
\end{equation}
belongs to $H^{\alpha }$ (see \cite{tri},$\S$.13.1).\\

We use notations from definition \ref{scalaw}

\begin{thm}
Let $0<\alpha <1$ and $f\in H^{\alpha } (c,C)$ defined on an open interval $U\subset \rR$
such that $I=[0,1] \subset U$. For all $\epsilon >0$, we define
\begin{equation}
l_{\epsilon } =\epsilon^{\alpha -1} \di \sqrt{\epsilon^{2(1-\alpha )} +c^2} ,\
L_{\epsilon } =\epsilon^{\alpha -1} \di \sqrt{\epsilon^{2(1-\alpha )} +C^2} .
\end{equation}
A scale law for $f$ is then given by
\begin{equation}
E(y,t)=(\alpha -1) (y -1/y ) .
\end{equation}
\end{thm}

{\it Proof}. We have
\begin{equation}
{\cal L}_{\epsilon} =\di\sum_{i=0}^{[1/\epsilon ]} \sqrt{\epsilon^2 + (f(\epsilon (i+1) )
-f(i\epsilon ))^2 }+
\di\sqrt{ (1-\epsilon [1/\epsilon ] )^2 +(f(1)-f(\epsilon [1/\epsilon ] ) )^2 } .
\end{equation}
The second term goes to zero when $\epsilon$ goes to zero. As $f \in H^{\alpha }$, we
have
\begin{equation}
c^2 \epsilon^{2\alpha } \leq
(f(\epsilon (i+1) )
-f(i\epsilon ))^2 \leq C^2 \epsilon^{2\alpha} .
\end{equation}
As a consequence, we obtain, for $\epsilon$ sufficiently small
\begin{equation}
\di\sum_{i=1}^{[1/\epsilon ]} \epsilon^{\alpha} \sqrt{\epsilon^{2(1-\alpha )} +c^2 }
+
\di\sqrt{ (1-\epsilon [1/\epsilon ] )^2 +(f(1)-f(\epsilon [1/\epsilon ] ) )^2 } .
\leq {\cal L}_{\epsilon} ,
\end{equation}
and
\begin{equation}
{\cal L}_{\epsilon } \leq
\di\sum_{i=1}^{[1/\epsilon ]} \epsilon^{\alpha} \sqrt{\epsilon^{2(1-\alpha )} +C^2 }
+\di\sqrt{ (1-\epsilon [1/\epsilon ] )^2 +(f(1)-f(\epsilon [1/\epsilon ] ) )^2 } .
\end{equation}
We deduce
\begin{equation}
\epsilon^{\alpha -1} \sqrt{\epsilon^{2(1-\alpha )} +c^2 }
\leq {\cal L}_{\epsilon }
\leq
\epsilon^{\alpha -1} \sqrt{\epsilon^{2(1-\alpha )} +C^2 }
.
\end{equation}
We then obtain
\begin{equation}
l_{\epsilon } \leq {\cal L}_{\epsilon } \leq L_{\epsilon } ,
\end{equation}
for $\epsilon$ sufficiently small. \\

By differentiating $L_{\epsilon}$ with respect to $\epsilon$, we obtain
\begin{equation}
\di {dL_{\epsilon } \over d\epsilon} =
\di {\alpha -1 \over \epsilon }
\left [
L_{\epsilon } -\di {1\over L_{\epsilon } }
\right ] .
\end{equation}
Using
\begin{equation}
\di {dL_{\epsilon } \over d\ln \epsilon} =\epsilon
\di {dL_{\epsilon } \over d\epsilon} ,
\end{equation}
we obtain the scale function $E(y,t)=(\alpha -1)(y-1/y)$. We verify that
$l_{\epsilon }$ satisfies the same differential equations. $\Box$ \\

The previous result is best analyzed in term of the new variables
\begin{equation}
x_{\epsilon} =1/l_{\epsilon } ,\ {\cal X}_{\epsilon } =1/{\cal L}_{\epsilon } ,\
X_{\epsilon } =1/L_{\epsilon } .
\end{equation}
When $\epsilon$ goes to $0$, we have $x_{\epsilon}$, $X_{\epsilon}$ and ${\cal X}_{\epsilon }$
which go to $0$. Moreover, the scale law for these new functions is
\begin{equation}
\label{slm}
\di {dx \over dt} =(1-\alpha ) (x-x^3 ) .
\end{equation}
Indeed, by making the change of variables $x=1/y$ in the scale law
\begin{equation}
\di {dy\over dt} =(\alpha -1) (y-1/y) ,
\end{equation}
and using the relation
\begin{equation}
\di {dx\over dt} = -(1/y^2 ) \di {dy\over dt} ,
\end{equation}
we obtain the result.\\

The classical linearization theorem of Poincare (\cite{ar},\cite{hk}) allows us to find, in a
neighborhood of $x=0$, an analytic change of variables $z=h(x)$, such that
the differential equation (\ref{slm}) is transformed into
\begin{equation}
\label{lsl}
\di {dz\over dt} =(1-\alpha )  z .
\end{equation}
The set of H\"olderian functions $H^{\alpha }$ induces, up to analytic changes of variables,
the linear scale law (\ref{lsl}).

\begin{rema}
In Galilean scale relativity, the set of functions which admit a linear scale law allows us
to define the Djinn variable (see \cite{c2},$\S$.3.2.5 and $\S$.3.3).
\end{rema}

\subsubsection{Non uniform H\"olderian functions}

In this section, we consider non uniform H\"oderian functions.

\begin{defi}
Let $\alpha (t) :\rR \rightarrow ]0,1[$. We denote by $H^{\alpha (.)}$ the set
of continuous H\"olderian functions satisfying, for all $h>0$ sufficiently small,
\begin{equation}
c h^{\alpha (t)} < \mid f(t+h) -f(t) \mid < C h^{\alpha (t)} ,
\end{equation}
where $c>0$, $C>0$ are constants.
\end{defi}

For $H^{\alpha (.)}$ functions, we have not been able to derive a scale law. We then
introduce a weak notion :

\begin{defi}
\label{wsl}
We say that $f$ admits weak-scale laws, if there exists $l_{\epsilon}$ and $L_{\epsilon}$
such that
\begin{equation}
l_{\epsilon} \leq {\cal L}_{\epsilon} \leq L_{\epsilon},
\end{equation}
satisfying
\begin{equation}
{\cal L}_{\epsilon} =O(l_{\epsilon} ),\ \ {\cal L}_{\epsilon} =O( L_{\epsilon} ) ,
\end{equation}
and for which we can find two functions $E_- :\rR \times \rR \rightarrow \rR$ and
$E_+ :\rR \times \rR \rightarrow \rR$ such that
\begin{equation}
\di {dl_{\epsilon} \over d\ln \epsilon} =E_- (l_{\epsilon} ,\ln \epsilon ),\ \
\di {dL_{\epsilon} \over d\ln \epsilon} =E_+ (L_{\epsilon} ,\ln \epsilon ).
\end{equation}
\end{defi}

For all $\epsilon >0$, we define
\begin{equation}
\alpha_i (\epsilon )=\alpha (i/\epsilon i+1/\epsilon ) ,
\end{equation}
with $i=0,\dots ,[1/\epsilon ]$, where $[x]$ denotes the integer part of $x$. \\

We define
\begin{equation}
\label{expo}
\gamma (\epsilon )= \min_{i=0,\dots , [1/\epsilon ]} \alpha_i (\epsilon ) , \
\beta (\epsilon )=\max_{i=0,\dots , [1/\epsilon ]} \alpha_i (\epsilon ) .
\end{equation}

\begin{thm}
Let $f \in H^{\alpha (.)}$. For all $\epsilon >0$ sufficiently small we assume that
the exponents (\ref{expo}) are differentiable functions with respect to $\epsilon$.
Then, $f$ admits for weak-scale laws
\begin{equation}
E_- (x,t)=(1-\alpha -t \di\gamma ') (x-x^3),\ \
E_+ (x,t)=(1-\alpha -t \di \beta ') (x-x^3 ) .
\end{equation}
\end{thm}

\section{Scale calculus}

In \cite{bc1}, we introduce the notion of {\it quantum derivatives}. In this
section, we give a less rigid definition, which allow us to discuss
more easily algebraic properties of these operators.

\subsection{Left and right quantum difference operators}

Let $h>0$ be given. If $f$ possesses a non zero $h$-minimal resolution, then for
$\epsilon <\epsilon (f,h)$, one must take into account the non differentiable
character of $f$ with respect to its forward and backward mean representations. A
possible way to do it, is to say that the backward and forward derivatives of $f$
at $\epsilon$ carry different informations on the local behaviour of $f$. The
idea of quantum derivatives formalize this idea.

\begin{defi}
Let $h>0$, and $f$ be a continuous, real valued fonction.\\

If $\epsilon (f,h) >0$, for all $\epsilon (f,h) >\epsilon >0$, we define left and right
quantum difference operators of $f$ at point $t$, the
quantities
\begin{equation}
\td_+^{\epsilon} f (t) = \di {f(t+\epsilon )-f(t) \over \epsilon} ,\ \td_-^{\epsilon} f(t) =
\di {f(t)-f(t-\epsilon )\over \epsilon} ,
\end{equation}
respectively.\\

If $\epsilon (f,h) =0$, then
\begin{equation}
\Delta_+^0 f(t) =\Delta_-^0 f(t) =f' (t) .
\end{equation}
\end{defi}

{\bf Remark}. We can give a more rigid definition of quantum derivatives by fixing the
$\epsilon$ to be the minimal resolution of the function considered for a given
$h>0$.

\subsection{The scale difference operator}

The scale difference operator, first introduced in \cite{bc1} following Nottale's work \cite{no1},
is intended to summarize the information given by quantum difference operator, need to perform
the local analysis of a given non differentiable function.

\begin{defi}
Let $h>0$, and $f$ be a continuous function such that $\epsilon (f,h) \geq \epsilon >0$.
The $\epsilon$-scale difference operator
of $f$ at point $t$, is a complex operator, denoted by ${\Box_{\epsilon} /\Box t}$, defined by
\begin{equation}
\di {\Box_{\epsilon} f \over \Box t} (t)=\di {1\over 2} (\td_+^{\epsilon} f(t) +\td_-^{\epsilon} f(t))-
i \di {1\over 2} (\td_+^{\epsilon} f(t) -\td_-^{\epsilon} f(t)),\ i^2 =-1.
\end{equation}
\end{defi}

We have the following usefull property :

\begin{lem}[Gluing]
Let $f$ be a differentiable function. Then, we have
\begin{equation}
\di {\Box f\over \Box t} = \di {df\over dt} .
\end{equation}
\end{lem}

Moreover, if we denote by $f_{\epsilon} (t)$ the mean function
$f_{\epsilon} (t) = (1/2\epsilon ) \di\int_{t-\epsilon}^{t+\epsilon} f(s) ds$, we have that
\begin{equation}
\mbox{Re}(\di {\Box_{\epsilon} f\over \Box t} (t))=(f_{\epsilon} )' (t) ,
\end{equation}
where $Re$ denote the real part of a complex number. The imaginary part of $\Box_{\epsilon} f/\Box t$
is the fluctuation of the forward mean function with respect to the backward mean function.

\subsection{Main result}

We denote by
$H^{\alpha}$ the set of real valued functions which are H\"older and inverse H\"older of
exponent $\alpha$, which mean that for all $\epsilon$ sufficiently small, and
$\mid t' -t \mid \leq \epsilon$, there exists $c>0$ and $C>0$, such that
\begin{equation}
c\epsilon^{\alpha} \leq \mid f(t')-f(t)\mid \leq C \epsilon^{\alpha} .
\end{equation}
An important theorem, from the point of view of the scale relativity, is the following :

\begin{thm}
Let $h>0$, $f(x,t)$ be a $C^n$ function, $x(t)$ a continuous function such that
$\epsilon (x,h) >0$. Then, for $\epsilon (x,h) \geq \epsilon >0$, $\epsilon$ sufficiently
small, we have
\begin{equation}
\di {\Box_{\epsilon} f (x(t),t) \over \Box t} =
\di {\partial f\over \partial t} +
\di \sum_{j=1}^n \di {1\over j!} \di {\partial^j f \over \partial x^j } (x(t),t)
\epsilon^{j-1} a_{\epsilon ,j} (t) +o (\epsilon^{1/n} ),
\end{equation}
where
\begin{equation}
a_{\epsilon ,j} (t) =
\di {1\over 2}
\left [
\left (
(\Delta_+^{\epsilon} x )^j -(-1)^j (\Delta_-^{\epsilon} x )^j
\right )
-i
\left (
(\Delta_+^{\epsilon} x )^j +(-1)^j (\Delta_-^{\epsilon} x )^j
\right )
\right ]
.
\end{equation}
\end{thm}

The proof follows from the following lemma :

\begin{lem}
Let $f(x,t)$ be a real valued function of class $C^{n+1}$, $n\geq 3$, and let $X(t)$ be a
continuous real valued function of class $H^{1/n}$. For $\epsilon$ sufficiently small,
the right and left derivatives of $f(X(t),t)$ are given by
\begin{equation}
\td_{\sigma}^{\epsilon} f(X(t),t) =\di {\partial f \over \partial t} (X(t),t) +\sigma \di\sum_{i=1}^n
\di {1\over i!} \di {\partial^i f \over \partial x^i} (X(t),t) \epsilon^{-1} (\sigma \epsilon \td_{\sigma}^{\epsilon} X(t))^i
+o(\epsilon^{1/n} ) ,
\end{equation}
for $\sigma =\pm$.
\end{lem}

For $n=2$, we obtain the so-called It\^o formula :
\begin{equation}
\left .
\begin{array}{lll}
\td_{\sigma}^{\epsilon} f(X(t),t) & = & \di {\partial f \over \partial t} (X(t),t)
+\di {\partial f \over \partial x} (X(t),t)
\td_{\sigma}^{\epsilon} X(t) \\
 & & +\di {1\over 2} {\partial^2 f \over \partial x^2} (X(t),t)
\epsilon (\td_{\sigma}^{\epsilon} X(t))^2
+o(\epsilon^{1/2} ) ,
\end{array}
\right .
\end{equation}

{\it Proof}. This follows from easy computations. First, we remark that, as $X(t)\in H^{1/n}$,
we have $\mid \epsilon \td_{\sigma}^{\epsilon} X(t) \mid =o(\epsilon^{1/n} )$. Moreover,
\begin{equation}
f(X(t+\epsilon ) ,t+\epsilon ) =f(X(t)+\epsilon \td_+^{\epsilon} X(t), t+\epsilon ) .
\end{equation}
By the previous remark, and the fact that $f$ is of order $C^{n+1}$, we can make a Taylor
expansion up to order $n$ with a controled remainder.
\begin{equation}
\left .
\begin{array}{lll}
f(X(t+\epsilon) ,t+\epsilon )& = & f(X(t),t)+\di\sum_{k=1}^{n} \di {1\over k!} \di \sum_{i+j=k}
(\epsilon \td_+^{\epsilon} X(t))^i \epsilon^j \di {\partial^k f \over \partial^i x \partial^j t} (X(t),t) \\
 & & +o((\epsilon \td_+^{\epsilon} X(t))^{n+1} ) .
\end{array}
\right .
\end{equation}
As a consequence, we have
\begin{equation}
\epsilon \td_+^{\epsilon} f(X(t),t) =\di\sum_{k=1}^{n} \di {1\over k!} \di \sum_{i+j=k}
(\epsilon \td_+^{\epsilon} X(t))^i \epsilon^j \di {\partial^k f \over \partial^i x \partial^j t} (X(t),t)
+o((\epsilon \td_+^{\epsilon} X(t))^{n+1} ) .
\end{equation}
By selecting terms of order less or equal to one in $\epsilon$ in the right of this equation, we
obtain
\begin{equation}
\left .
\begin{array}{lll}
\epsilon \td_+^{\epsilon} f(X(t),t) & = & \epsilon \left [
\di {\partial f \over \partial t} (X(t),t) +\di\sum_{i=1}^n
\di {1\over i!} \di {\partial^i f \over \partial x^i} (X(t),t) \epsilon^{-1} (\epsilon \td_+^{\epsilon} X(t))^i
\right ] \\
 & & +o(\epsilon^2 \td_+^{\epsilon} X(t) ) .
\end{array}
\right .
\end{equation}
Dividing by $\epsilon$, we obtain the lemma. $\Box$

\subsection{The complex case}

In part II, we need to apply the scale operator to {\it complex valued} functions. We
extend the definition of $\Box_{\epsilon } /\Box t$ in order to cover this case.\\

In the following, if $z\in \cC$, we denote by $\re (z)$ and $\im (z)$, the real and imaginary
part of $z$.

\begin{defi}
Let $h>0$ and ${\cal C} (t)$ be a complex valued function. We denote ${\cal C} (t) ={\cal C}_r (t) +i
{\cal C}_m (t)$, where ${\cal C}_r (t)=\re {\cal C}(t)$ and ${\cal C}_m (t)=\im {\cal C}(t)$.
We define $\Box_{\epsilon} {\cal C} / \Box t$ by
\begin{equation}
\di {\Box_{\epsilon} {\cal C} \over \Box t} =\di {\Box_{\epsilon} {\cal C}_r \over \Box t}
+i \di {\Box_{\epsilon} {\cal C}_m \over \Box t} ,
\end{equation}
for $0<\epsilon <\min (\epsilon ({\cal C}_r ,h),\epsilon ({\cal C}_m ,h))$.
\end{defi}

We then have :

\begin{lem}
Let $h>0$ and $C(x,t) :\rR \times \rR \rightarrow \cC$ be a $C^n$ complex
valued function. Let $x(t)$ be a continuous function such that $\epsilon (x,h) >0$.
We denote ${\cal C}(t)=C(x(t),t)$. Then,
for $0<\epsilon \leq \epsilon (x,h)$ sufficiently small, we have
\begin{equation}
\di {\Box_{\epsilon} {\cal C} \over \Box t} =
\di {\partial C\over \partial t} +\di {\Box_{\epsilon} x \over \Box t}
\di {\partial C\over \partial x} +\di \sum_{j=2}^n {1\over j!} a_{\epsilon ,j} (t) \di
{\partial^j C \over \partial x^j} \epsilon^{j-1} +o(\epsilon^{1/n} ),
\end{equation}
where
\begin{equation}
a_{\epsilon ,j} (t) =\di {1\over 2} [(\Delta_+^{\epsilon} x)^j -(-1)^j (\Delta_-^{\epsilon} x
)^j ] -i \di {1\over 2}
[(\Delta_+^{\epsilon} x)^j +(-1)^j (\Delta_-^{\epsilon} x )^j ] , \ j=2,\dots ,n .
\end{equation}
\end{lem}

\section{Algebraic properties of quantum difference operators}

In the following, we denote by $C^0$ the set of continuous real valued functions.

\begin{lem}
For all $\epsilon >0$, and all $f\in C^0$ and $g\in C^0$, we have : \\

i) $\td_{\sigma}^{\epsilon} (f+g)=\td_{\sigma} f +\td_{\sigma} g$, $\sigma =\pm$, \\

ii) For all $\lambda \in \rR$, $\td (\lambda f )=\lambda \td_{\sigma} f$, $\sigma =\pm$.
\end{lem}

The proof is straightforward and let to the reader. \\

Our main is to compare quantum derivatives to classical derivatives. The main
properties of classical derivatives is the so called {\it Leibniz rule}, which says
that, $(fg)'=f'g+fg'$. In our case, we have a more complicated formula.

\begin{lem}
\label{r1}
For all $\epsilon >0$, and all $f\in C^0$, $g\in C^0$, we have for $\sigma =\pm$,
\begin{equation}
\label{res}
\td_{\sigma}^{\epsilon} (fg)(x)  =\td_{\sigma}^{\epsilon} f(x) g(x)+f(x) \td_{\sigma}^{\epsilon} g(x) +\sigma \epsilon
\td_{\sigma}^{\epsilon} f(x) \td_{\sigma}^{\epsilon} g(x) .
\end{equation}
\end{lem}

{\it Proof}. Easy computations lead to the following formulas :
\begin{equation}
\label{r0}
\left .
\begin{array}{lll}
\td_+^{\epsilon} (fg)(x) & = & \td_+^{\epsilon} f(x) g(x+\epsilon )+f(x) \td_+^{\epsilon} g(x) , \\
\td_-^{\epsilon} (fg)(x) & = & \td_-^{\epsilon} f(x) g(x) +f(x-\epsilon ) \td_-^{\epsilon} g(x) .
\end{array}
\right .
\end{equation}
By definition of the quantum derivatives, we have
\begin{equation}
f(x+\epsilon )=f(x)+\epsilon \td_+^{\epsilon} f(x),\ \ g(x-\epsilon )=g(x)-\epsilon \td_-^{\epsilon} g(x) .
\end{equation}
By replacing $f(x+\epsilon )$ and $g(x-\epsilon )$ in (\ref{r0}), we obtain (\ref{res}). $\Box$ \\

Of course, one can derived others formulas. However, equation (\ref{res}) is the most symetric one.
Here, we give another expression. \\

For all $\epsilon$, we denote by $\tau_{\epsilon} :C^0 \rightarrow C^0$ the classical translation, defined by
\begin{equation}
\tau_{\epsilon} (f) (x)=f(x+\epsilon ) ,\ \forall x.
\end{equation}
We have the following lemma :

\begin{lem}
For all $\epsilon$, we have $\td_-^{\epsilon} \circ \tau_{\epsilon} =\td_+^{\epsilon}$ and
$\tau_{\epsilon} \circ \td_{\sigma}^{\epsilon} =\td_{\sigma}^{\epsilon} \circ \tau_{\epsilon }$.
\end{lem}

As a consequence, we obtain the following version of lemma \ref{r1} :

\begin{lem}
\label{r2}
For all $\epsilon >0$, and all $f\in C^0$, $g\in C^0$, we have
\begin{equation}
\left .
\begin{array}{lll}
\td_+^{\epsilon} (fg)(x) & = & \td_+^{\epsilon} f(x) h_{\epsilon} (x)+f(x) \td_-^{\epsilon} h_{\epsilon} (x) , \\
\td_-^{\epsilon} (fg)(x) & = & \td_-^{\epsilon} f(x) v_{\epsilon} (x) +f(x) \td_+^{\epsilon} v_{\epsilon} (x) ,
\end{array}
\right .
\end{equation}
where $h_{\epsilon} =\tau_{\epsilon} \circ g$ and $v_{\epsilon} =\tau_{-\epsilon} \circ g$.
\end{lem}

In the following, we discuss what are the fundamental differences between quantum derivatives and
classical derivations.

\subsection{Derivations}

We recall that an operator $D$ on an abstract algebra $(A,.)$, is a derivation if for all $(f,g)\in A^2$,
it satisfies the Leibniz relation $D(fg)=Df.g+f.Dg$. We refer to Jacobson (\cite{ja}, chap.1, $\S$.2, p.7-8)
for more details. \\

We denote by $Der(A)$ the set of derivations on $A$. $Der(A)$ is a vector space, but not an algebra. However,
by posing $[D_1 ,D_2]=D_1 D_2 -D_2 D_1$, the usual Lie Bracket, the set $(Der(A), [.,.])$ is a Lie
algebra. \\

We denote by $\rR \langle\langle Der(A)\rangle\rangle$ the ring of formal power series on the alphabet
$Der(A)$. We can defined a coalgebra structure on $\rR \langle\langle Der(A)\rangle\rangle$. We refer to
Bourbaki (\cite{bo},chap.3) for more details about coalgebras and bialgebras. \\

Let $u: \rR \langle\langle Der(A)\rangle\rangle \rightarrow \rR$ be the homomorphism associating to each
serie its constant term. \\

For each $D\in Der(A)$ we define a linear map $\Delta : Der(A) \rightarrow Der(A)\otimes Der(A)$, by
$\Delta (D)=D\otimes 1 +1\otimes D$. Then, the following diagram commutes
\begin{equation}
\left .
\begin{array}{lll}
A\otimes A & \stackrel{\Delta (D)}{\rightarrow } & A\otimes A , \\
\nu \downarrow & & \downarrow \nu \\
A & \stackrel{D}{\rightarrow } & A ,
\end{array}
\right .
\end{equation}
where $\nu$ is the natural morphism defined by
\begin{equation}
\nu :
\left .
\begin{array}{lll}
A\otimes A & \rightarrow A , \\
f\otimes g & \mapsto & f. g.
\end{array}
\right .
\end{equation}
We can extend $\Delta$ such that for each $\Delta (D_1 D_2 )=\Delta (D_1 )\Delta (D_2 )$, $D_1$, $D_2 \in
Der (A)$, and the usual product on $\rR \langle\langle Der(A)\rangle\rangle
\otimes \rR \langle\langle Der(A)\rangle\rangle$, $(a\otimes b).(c\otimes d)=(ac\otimes bd)$.
For each word, $D_1 \dots D_r$, we define $\Delta (D_1 \dots D_r)=\Delta(D_1 )\Delta (D_2 \dots D_r)$.
We then extend $\Delta$ to $\rR \langle\langle Der(A)\rangle\rangle$ by linearity. With this
definition of $\Delta$, for each $S\in \rR \langle\langle Der(A)\rangle\rangle$, the following diagram commutes
\begin{equation}
\left .
\begin{array}{lll}
A\otimes A & \stackrel{\Delta (S)}{\rightarrow } & A\otimes A , \\
\nu \downarrow & & \downarrow \nu \\
A & \stackrel{S}{\rightarrow } & A ,
\end{array}
\right .
\end{equation}
As a consequence, we have the following lemma :

\begin{lem}
The triple $(\rR \langle\langle Der(A)\rangle\rangle ,\Delta ,u)$ is a bialgebra.
\end{lem}

\subsection{Quantum bialgebra}

We follow the previous section on derivations.

\begin{defi}
For all $\epsilon >0$, and $\sigma =\pm$, we denote by $\td_{\sigma}^{\epsilon}$ an operator
acting on $C^0 \otimes C^0$, where $\otimes$ is the classical tensor product, and defined by
\begin{equation}
\td_{\sigma}^{\epsilon} :
\left .
\begin{array}{lll}
C^0 \otimes C^0 & \rightarrow & C^0 \otimes C^0 , \\
f\otimes g & \mapsto & \td_{\sigma}^{\epsilon} f \otimes g +
f\otimes \td_{\sigma}^{\epsilon} g +\sigma \epsilon \td_{\sigma}^{\epsilon} f \otimes
\td_{\sigma}^{\epsilon} g.
\end{array}
\right .
\end{equation}
\end{defi}

Let $\Omega_{\epsilon} =\{ \td_+^{\epsilon} ,\td_-^{\epsilon} \}$ be the alphabet of two letters
$\td_+^{\epsilon}$ and $\td_-^{\epsilon} \}$. We denote by $\Omega_{\epsilon}^*$ the set of words
$\un{\omega} =\omega_1 \dots \omega_n$, $\omega_i \in \Omega_{\epsilon }$ for all $i=1,\dots ,n$, where
$\omega_i \omega_j$ denote the natural composition of operators. For example, a possible word is
$\td_+^{\epsilon} \td_-^{\epsilon}$. \\

{\bf Remark}. Here we consider the alphabet of quantum operators for a fixed $\epsilon >0$. It is
possible that some particular problems of scale relativity require a complete alphabet
$\Omega =\{ \td_{\sigma}^{\epsilon},\sigma =\pm ,\epsilon >0 \}$. \\

We denote by $A_{\epsilon} =\rR \langle \langle \Omega_{\epsilon} \rangle \rangle$ the algebra of
formal power series constructed on $\Omega_{\epsilon}$ (with its classical algebraic structure).\\

We can define a linear map from $A_{\epsilon }$ to $A_{\epsilon} \times_{\rR} A_{\epsilon}$, denoted
$\Delta$. First, we define $\Delta$ on $\td_{\sigma}^{\epsilon}$ by
\begin{equation}
\left .
\begin{array}{lll}
A_{\epsilon} & \stackrel{\Delta}{\rightarrow} & A_{\epsilon}\otimes A_{\epsilon}, \\
\td_{\sigma}^{\epsilon} & \mapsto & \td_{\sigma}^{\epsilon} \otimes I +I\otimes \td_{\sigma}^{\epsilon}
+\sigma \epsilon \td_{\sigma}^{\epsilon} \otimes \td_{\sigma}^{\epsilon} .
\end{array}
\right .
\end{equation}
We have the following equality
\begin{equation}
\td_{\sigma}^{\epsilon} \circ \nu = \nu \circ \Delta (\td_{\sigma}^{\epsilon} ),
\end{equation}
which is equivalent to the commutativity of the diagramm
\begin{equation}
\left .
\begin{array}{lll}
C^0 \otimes C^0 & \stackrel{\Delta (\td_{\sigma}^{\epsilon} )}{\rightarrow} & C^0 \otimes C^0 , \\
\downarrow \nu & & \downarrow \nu  \\
C^0 & \stackrel{\td_{\sigma}^{\epsilon}}{\rightarrow} & C^0 .
\end{array}
\right .
\end{equation}
We also define $\Delta (I)=I\otimes I$ in order that $I \circ \nu =\nu \circ \Delta (I)$. \\

We extend $\Delta$ to $A_{\epsilon}$ by linearity.

\begin{lem}
The linear map $\Delta$ is an algebra homomorphism.
\end{lem}

{\it Proof}. The proof is done by induction. Let $\td_{\sigma}^{\epsilon}$ and
$\td_{\sigma'}^{\epsilon}$ be two letters of $\Omega_{\epsilon}$. We have
\begin{equation}
\left .
\begin{array}{lll}
\td_{\sigma}^{\epsilon} \td_{\sigma'}^{\epsilon} (fg) & = &
\td_{\sigma}^{\epsilon} \td_{\sigma'}^{\epsilon} f .g +f. \td_{\sigma}^{\epsilon} \td_{\sigma'}^{\epsilon} g \\
 & & + \td_{\sigma}^{\epsilon} f .\td_{\sigma'}^{\epsilon} g +\td_{\sigma'}^{\epsilon} f .\td_{\sigma}^{\epsilon} g \\
 & & +\epsilon [ \td_{\sigma}^{\epsilon} \td_{\sigma'}^{\epsilon} f .\td_{\sigma'}^{\epsilon} g
 +\td_{\sigma'}^{\epsilon} f . \td_{\sigma}^{\epsilon} \td_{\sigma'}^{\epsilon} g \\
 & & +\td_{\sigma}^{\epsilon} \td_{\sigma'}^{\epsilon} f .\td_{\sigma}^{\epsilon} g
  +\td_{\sigma}^{\epsilon} f . \td_{\sigma}^{\epsilon} \td_{\sigma'}^{\epsilon} g ] \\
 & & +\epsilon^2 \td_{\sigma}^{\epsilon} \td_{\sigma'}^{\epsilon} f . \td_{\sigma}^{\epsilon} \td_{\sigma'}^{\epsilon}
g .
\end{array}
\right .
\end{equation}
As a consequence, we have
\begin{equation}
\left .
\begin{array}{lll}
\Delta (\td_{\sigma}^{\epsilon} \td_{\sigma'}^{\epsilon} ) & = &
\td_{\sigma}^{\epsilon} \td_{\sigma'}^{\epsilon} \otimes I +I\otimes  \td_{\sigma}^{\epsilon} \td_{\sigma'}^{\epsilon}  \\
 & & + \td_{\sigma}^{\epsilon} \otimes \td_{\sigma'}^{\epsilon} +\td_{\sigma'}^{\epsilon} \otimes \td_{\sigma}^{\epsilon} \\
 & & +\epsilon [ \td_{\sigma}^{\epsilon} \td_{\sigma'}^{\epsilon} \otimes \td_{\sigma'}^{\epsilon}
 +\td_{\sigma'}^{\epsilon} \otimes \td_{\sigma}^{\epsilon} \td_{\sigma'}^{\epsilon}  \\
 & & +\td_{\sigma}^{\epsilon} \td_{\sigma'}^{\epsilon} \otimes \td_{\sigma}^{\epsilon}
  +\td_{\sigma}^{\epsilon} \otimes \td_{\sigma}^{\epsilon} \td_{\sigma'}^{\epsilon}  ] \\
 & & +\epsilon^2 \td_{\sigma}^{\epsilon} \td_{\sigma'}^{\epsilon} \otimes \td_{\sigma}^{\epsilon} \td_{\sigma'}^{\epsilon} .
\end{array}
\right .
\end{equation}
By definition, the classical product on $A_{\epsilon} \otimes A_{\epsilon}$ is defined by : for all $a\otimes b$ and
$c\otimes d$ in $A_{\epsilon}\otimes A_{\epsilon}$, we have
\begin{equation}
(a\otimes b). (c\otimes d)=ac \otimes bd .
\end{equation}
An easy computation proves that
\begin{equation}
\Delta (\td_{\sigma}^{\epsilon} \td_{\sigma'}^{\epsilon} )=\Delta (\td_{\sigma}^{\epsilon} )
\Delta (\td_{\sigma'}^{\epsilon} ) .
\end{equation}
By induction and linearity, we obtain the lemma. $\Box$ \\

We define an algebra homomorphism from $A_{\epsilon }$ to $\rR$, denoted by $u$, by associating to each
formal power series its constant term. \\

With $u : A_{\epsilon } \rightarrow \rR$ and $\Delta : A_{\epsilon} \rightarrow A_{\epsilon}
\otimes A_{\epsilon}$, we define a coalgebra structure on $A_{\epsilon}$.

\begin{lem}
The triple $(A_{\epsilon}, \Delta , u )$ is a coalgebra.
\end{lem}

Moreover, as $u$ and $\Delta$ are homomorphism, we obtain the stronger result :

\begin{lem}
The triple $(A_{\epsilon}, \Delta , u )$ is a bialgebra.
\end{lem}

We have not found a natural graduation on this bialgebra.\\

{\bf Remark}. 1. It will be interseting to discuss the possible relation to {\it Quantum groups} introduced by
Drinfeld (\cite{dr1},\cite{dr2}) and quasi-triangular Hopf algebra. We introduce the
natural commutativity involution
\begin{equation}
\tau :
\left .
\begin{array}{lll}
A_{\epsilon} \otimes A_{\epsilon} & \rightarrow & A_{\epsilon } \otimes A_{\epsilon} , \\
 a\otimes b & \longmapsto & b\otimes a .
\end{array}
\right .
\end{equation}
Let $\Delta^{op} =\tau \circ \Delta$. The map $\tau$ is an algebra automorphism of
$A_{\epsilon} \otimes A_{\epsilon }$ and the following diagram commutes
\begin{equation}
\left .
\begin{array}{lll}
A_{\epsilon} & \stackrel{\Delta}{\longrightarrow} & A_{\epsilon} \otimes A_{\epsilon} , \\
id_{A_{\epsilon}} \downarrow &  & \downarrow \tau \\
A_{\epsilon} & \stackrel{\Delta^{op}}{\longrightarrow} & A_{\epsilon } \otimes A_{\epsilon} .
\end{array}
\right .
\end{equation}
If $(A_{\epsilon } ,\Delta ,u)$ is a quasitriangular algebra, then, following \cite{ra}, there
exists an invertible element of $A_{\epsilon} \otimes A_{\epsilon}$ such that
\begin{equation}
\Delta^{op} (a) =R \Delta (a) R^{-1},\ \mbox{for all}\ a\in A_{\epsilon} .
\end{equation}
In our case, we easily have $\Delta^{op} =\Delta$ so that $(A_{\epsilon} ,\Delta ,u)$ so that
$R$ is trivial. \\

2. A natural idea is to consider the quantum bialgebra as a ``deformation" of the classical
Hopf algebra associated to derivations. The word derivation must be taken with care, because
a notion of deformation for Hopf algebra already exists and it is not clear if this is
the good one to consider here. \\

3. There is no natural extension of the quantum bialgebra in order to take into account the
scale derivative. The basic problem being that the scale derivative is a complex valued operator on
real valued functions. Then, we have no natural composition of these operators.

\subsection{A remark on Rieman-Liouville fractional calculus and the local fractional
calculus}

A basic way to deal with non differentiable functions is to use {\it fractional calculus}.
As an example, one can consider the classical (left and right) {\it Riemann-Liouville derivative}, defined
by
\begin{equation}
D_z^{\alpha} f(x) =\di {\Gamma (\alpha +1)\over 2\pi i} \di \int_0^{z^+}
f(t) (t-z)^{-\alpha -1} dt ,
\end{equation}
for $\alpha \in \cC \setminus \{ -1 ,-2,\dots \}$. The Riemann-Liouville derivative
is not at all a {\it derivation} on the algebra of continuous functions. Indeed, one
has (see \cite{Os}) :
\begin{equation}
D_z^{\alpha} (uv)=\di\sum_{n=0}^{\infty} \left (
\begin{array}{c}
\alpha \\
n
\end{array}
\right )
D_z^{\alpha -n} u D_z^n v .
\end{equation}
Moreover, this operator, which is a direct analytic generalization of Cauchy formula
doesn't have a clear geometrical interpretation, despite recent advances (see \cite{ben}). \\

In \cite{bc2}, we have obtain, following a previous work of Kolvankar and
Gangal \cite{kolga}, a derivation, called the (left and right) {\it local fractional
derivative}, by {\it localyzing} the (left and right) Riemann-Liouville derivative.\\

We have prove that this localization takes a simple form.

\begin{defi}
The right and left local fractional derivative of $f$ at point $x_0$ of order $\alpha$
is defined by
\begin{equation}
\lim_{x\rightarrow x_0^+} \di {f(x)-f(x_0 )\over (x-x_0)^{\alpha}} ,\ \mbox{\rm and}\ \
\lim_{x\rightarrow x_0^-} \di {f(x_0)-f(x) \over (x_0 -x)^{\alpha}} ,
\end{equation}
respectively.
\end{defi}

We have introduce the fractional derivative of $f$ at point $x_0$ by collecting the
two quantities $d^{\alpha}_+ f(x_0)$ and $d^{\alpha}_- f(x_0 )$ in a single
quantity, i.e.
\begin{equation}
d^{\alpha} f(x_0 ) =\di {1\over 2} \left (
d^{\alpha}_+ f(x_0 ) + d^{\alpha}_- f(x_0 ) \right )
+i \di {1\over 2}
\left (
d^{\alpha}_+ f(x_0 ) + d^{\alpha}_- f(x_0 ) \right ) .
\end{equation}

Moreover, we obtain a clear geometrical meaning by connecting the exponent $\alpha$ of
differentiation to the local maximal H\"older regularity of the curve.\\

However, such a derivative has important problems. First, there exist no {\it integral operator}
contrary to the Riemann-Liouville fractional derivative where there exists the Riemann-Liouville
integral. Moreover, the set of points on which the local fractional derivative is non zero is in
most of the case trivial, i.e. of zero measure. Precisely, we have the following
theorem.

\begin{thm}
\label{conjec}
The fractional differential equations of the form $d^{\alpha} f(x)=a(x)+ib(x)$, $0<\alpha <1$, where
$a(x)$ and $b(x)$ are continuous functions such that there exists $x_0 \in \rR$ such that
$\mid a(x_0) \mid \not= \mid b(x_0 )\mid$, have no solutions.
\end{thm}

This theorem solves a conjecture of \cite{bc4}. It must be pointed out that the condition
$\mid a(x_0)\mid \not= \mid b(x_0 )\mid$ is generic.\\

{\it Proof}. We have $d^{\alpha}_+ f(x)=c(x)$ and $d^{\alpha}_- f(x)=d(x)$, where
$c(x)=a(x)+b(x)$ and $d(x)=a(x)-b(x)$ are continuous functions. By assumption, we have
$c(x_0 )\not =0$ and $d(x_0 )\not= 0$.

Let us assume that $c(x_0 )>0$ (the case where $c(x_0 )<0$ is similar).
By continuity, there exits $\epsilon >0$ such that
for all $x$ in the open interval $I_{\epsilon} (x_0 ) =]x_0 -\epsilon ,x_0 +\epsilon [$, we
have $d^{\alpha} f(x) >0$. Two case must be considered :\\

i) $d^{\alpha}_- f(x_0 )>0$ and ii) $d^{\alpha}_- f(x_0 )<0$. \\

In the case i), we define an open interval $J_{\epsilon} (x_0 )$ such that for all $x\in
J_{\epsilon} (x_0 )$ we have $d^{\alpha} f(x) >0$. As a consequence, the function $f$
is H\"olderian of exponent $\alpha$ (see \cite{bc1},theorem 3.9). Moreover, the function $f$
is injective on the interval $x\in K_{\epsilon} (x_0 )=J_{\epsilon} (x_0 )\bigcap I_{\epsilon} (x_0 )$. Indeed, if there
exists $x_1 ,x_2\in K_{\epsilon} (x_0 )$ such that $f(x_1 )=f(x_2 )$, then by
the fractional Rolle's theorem (\cite{bc1}), there $x\in K_{\epsilon} (x_0 )$ such that
$d^{\alpha}_+ f(x) d^{\alpha}_- f(x) \leq 0$, which is impossible by assumption. A
continuous function wich is injective is strictly monotone (see \cite{hawa},lemma 3.8,p.207).
But, a monotone function is almost everywhere differentiable by Lebesgue theorem (\cite{kf},p.319),
in contradiction with the assumption that $0<\alpha <1$.\\

In the case ii), we define an open interval $J_{\epsilon} (x_0 )$ such that for all
$x\in J_{\epsilon} (x_0 )$ we have $d^{\alpha} f(x) <0$. Again the function $f$ is an
H\"olderian function of exponent $\alpha$. Moreover, by (\cite{bc1},theorem 4.9),
all points in $K_{\epsilon} (x_0 )=J_{\epsilon} (x_0 )\bigcap I_{\epsilon} (x_0 )$ are
local minima. This is possible if and only if $f$ is a constant function, i.e. a
differentiable function in contradiction with the fact that $f$ is H\"olderian of exponent
$0<\alpha <1$. This concludes the proof of the theorem. $\Box$\\

We can characterize the {\it spectrum} of $f$, i.e. the set of values taken by the fractional
derivative of a continuous function.

\begin{cor}
The spectrum of a given $\alpha$-differentiable function is discontinuous or zero.
\end{cor}

{\it Proof}. The spectrum can't be continuous and non zero by theorem \ref{conjec}. As a
consequence, the spectrum can be zero or dicontinuous. $\Box$\\

We can be more precise on the nature of this spectrum.

\begin{thm}
\label{condis}
The set of fractional differential equations $d^{\alpha} f(x)=a(x)+ib(x)$, where
$0<\alpha <1$ and the functions $c(x)=a(x)+b(x)$, $d(x)=a(x)-b(x)$ keep a constant sign
on a given interval of $\rR$, has no solutions.
\end{thm}

{\it Proof}. This is the same proof as theorem \ref{conjec}. The continuity assumption of
theorem \ref{conjec} being only here to construct an interval where $c(x)$ and $d(x)$
keep a constant sign. $\Box$\\

As a consequence, the spectrum is discontinuous but of very special form, as we can't find
any interval of $\rR$ on which we have a constant sign. We define the {\it Dirichlet function}
as follows
\begin{equation}
{\cal D} (x) =\left \{
\begin{array}{cc}
1 & \ \ \mbox{\rm if}\ \ x\in \qQ , \\
-1 & \ \ \mbox{\rm if}\ \ x\in \rR \setminus \qQ .
\end{array}
\right .
\end{equation}
Then, the function $a(x)$ and $b(x)$ can be taken such that
\begin{equation}
a(x) =\di {3\over 2}{\cal D} (x) ,\ \ b(x)=-\di {1\over 2}{\cal D} (x) .
\end{equation}

We can extend all the previous theorem to the case where the order of fractional differentiation
is non constant, but a function of $x$.\\

In \cite{bc2}, we have derived , in the fractional calulus framework, the Schr\"odinger equation
from the Newton's equation of dynamics under the assumption that one-dimensional
quantum mechanical trajectories satisfy
\begin{equation}
\label{sccond}
(d^{\alpha}_{\sigma} f(x) )^2 =2\bar{h}/2m ,\ \sigma=\pm ,
\end{equation}
where
$\bar{h}$ is the reduced Planck constant $h/2\pi$ and $m$ is the mass of the particle. Of
course, this is impossible by theorem \ref{conjec}, as already proved in \cite{bc4}.
But as the Schr\"odinger equation is a well established equation of physics, we propose
in \cite{bc4} to consider a small perturbation of condition (\ref{sccond}), like
\begin{equation}
\label{second2}
(d^{\alpha}_{\sigma} f(x) )^2 =2\bar{h}/2m +\epsilon a_{\sigma} (t,f(t)),\ \sigma =\pm
\end{equation}
where $0<\epsilon <<1$, in order to permit the existence of non trivial solutions and to
obtain a small perturbation of the Schr\"odinger equation.\\

However, condition (\ref{second2}) lead to a deadlock. Indeed, for $\epsilon$ sufficiently
small, the quantities $d^{\alpha}_{\sigma} f(x)$ keep a
constant sign. As a consequence, by theorem \ref{condis} we have no solutions.\\

All these problems are solved in part II by using the scale calculus framework.

\section{Quantum representation of non differentiable functions}

We introduce the notion of quantum geometric representation for a continuous
non differentiable function. This notion is associated to minimal
resolution and the scale derivative. It turns out that a geometric
space displaying the basic features of the quantum geometric representation is
given by a simplified version of A. Connes formulation of the standard model of fundamental
interactions within the framework of noncommutative geometry.

\subsection{The quantum representation of a continuous non differentiable function}

Our previous results allow us to define a natural notion of scale derivative. The
scale derivative, which is a complex valued operator, contains the necessary informations
in order to perform a local analysis of a continuous non differentiable functions, and
take care of this non differentiability. \\

Let $f$ be a given continuous on differentiable function. The basic
functions associated to $f$, and from which we can deduce the scale derivative are the
forward and backward mean functions defined as $f_{\epsilon}^+ (t)=(1/2\epsilon) \di\int_t^{t+\epsilon} f(s) ds$ and
$f_{\epsilon}^- (t)=(1/2\epsilon ) \di\int_{t-\epsilon}^t f(s) ds$ respectively. \\

>From the geometrical view point, it means that in order to take into account the
non differentiable character of $f$, one must consider the disjoint union
$\Gamma_{\epsilon}^+ \cup \Gamma_{\epsilon}^-$,
where $\Gamma^{\sigma}$ is the graph of $f_{\epsilon}^{\sigma}$, $\sigma =\pm$.

\begin{defi}
Let $h>0$, $f$ be a continuous non differentiable function, and $\epsilon (f,h)$ be its
minimal resolution. For all $\epsilon >0$, the quantum geometric representation of $f$, denoted
$Q_{\epsilon} (f)$, is defined by : \\

i) For all $\epsilon >\epsilon (f,h)$, $Q_{\epsilon} (f)=\Gamma_{\epsilon }$, \\

ii) For all $0<\epsilon \leq \epsilon (f,h)$, $Q_{\epsilon} (f)=\Gamma_{\epsilon}^+ \cup \Gamma_{\epsilon}^-$, \\

where $\Gamma_{\epsilon}$, $\Gamma_{\epsilon}^+$, $\Gamma_{\epsilon}^-$, are the graph of
the mean function $f_{\epsilon}$, the forward mean function and the backward mean function
respectively.
\end{defi}

The non differentiability of $f$ induces a change in the geometric structure of the
geometric representation of $f$.\\

In the following, we consider graphs of real valued functions as submanifold of $\rR^2$.
As $f_{\epsilon}^+$ and $f_{\epsilon}^-$ are differentiable functions, one
deals with a disjoint union of differentiable submanifolds $\Gamma_{\epsilon}^+$ and $\Gamma_{\epsilon}^-$.
The basic features of the quantum representation of $f$, when $0<\epsilon < \epsilon (f,h)$, is that
$Q_{\epsilon}$ is composed by two differentiable submanifolds, which are close to each other, their
closeness being related to $\epsilon$. \\

As a consequence, a good understanding of the effects of a non differentiable function can be
obtained via the following simplified model : \\

Let $M$ be a one dimensional differentiable submanifold of $\rR^2$. Let $A=\{ a,b\}$ be a two points space.
We consider $Q=M\times A$. Then, $Q$ is the union of two copies of the manifold $M$ : $Q=M_a \cup M_b$. \\

{\bf Remark}: A more accurate model is the following : Let $A_{\epsilon}=\{ a_{\epsilon} , b_{\epsilon } \}$ be a two
points space, such that for $\epsilon >\epsilon (f,h)$, $A_{\epsilon}$ reduces to a point, i.e. $a_{\epsilon } =
b_{\epsilon}$. The simplified model is then $Q_{\epsilon } =M\times A_{\epsilon}$.

\subsection{Noncommutative geometry}

In this paragraph we only sketch a possible connexion between our point of view on
non differentiable functions and non commutative geometry. The idea is, by this way,
to obtain powerfull tools to study non differentiable functions which will be
relevant to physics. \\

In his book \cite{co}, A. Connes developp non commutative geometry. The basic idea is
to extend to the non commutative case the classical result of Gelfand and Naimark
relating $C^*$-algebras and locally compact spaces. \\

At the end of his book, Connes (\cite{co},p.568) discusses a particular example, where his theory already
lead to interesting results (this example is view as a preliminary step toward a
complete non commutative model for the standard model of quantum particles).\\

He considers a product of a differentiable manifold (the standard $\rR^4$) by a discrete
space, $A=\{ a,b\}$. Using non commutative geometry, he can make analysis on this
space. In particular, he defines a ``differential" operator which contains three
terms : the classical derivative on each copie of $\rR^4$ and a finite
difference. \\

It will be interesting to discuss the relevance of this construction with respect to
our approach to non differentiable functions.

\newpage
\part{Scale relativity and the Schr\"odinger equation}
\setcounter{section}{0}
\setcounter{defi}{0}
\setcounter{thm}{0}
\setcounter{cor}{0}
\setcounter{rema}{0}
\setcounter{lem}{0}
\setcounter{equation}{0}

\section{Introduction}

The Schr\"odinger equation is one of the basic pieces of quantum mechanics. Many
attempts already exist in order to derive it from the expected behaviours of
trajectories of quantum particles or from classical equations of the dynamics. We can cite
for example : \\

- Nelson stochastic approach \cite{nel}, \\

- Feynman perturbativ approach \cite{fey} ,\\

- Nottale's approach by the Scale relativity theory \cite{no1} .\\

In the following we discuss the derivation of the Schr\"odinger equation in the framework of
the Scale relativity of Nottale \cite{no1}. The main point is that,
contrary to Nelson or Feynman approach, it is based on a
{\it first principle}, namely, the {\it scale principle of relativity}, which is an extension of
the Einstein relativity principle to scales (of time and lengths).\\

The scale relativity principle introduced by Nottale's has a direct consequence on the
equations of the dynamics for a given particle. Indeed, they must keep the same form
under a scale transform, i.e. going from the classical scale to the atomic scale. Following
Feynman and Hibbs, the principal difference between the microscopic and macrocospic scale is
that typical paths become non differentiable. Then, we must be able to
transform the classical differential equations of the dynamics for functions which are
not at all differentiable.\\

This is done using the scale difference operator defined in part I. The scale relativity
principle is then equivalent to changing the classical derivative by the scale difference operator in the
Euler-Lagrange equations of the dynamics.\\

This quantification procedure called the scale quantization can be precisely defined
in $\S$.\ref{scalemap}, by introducing a quantization map, associating to each classical variables and
differential operators its
quantum counterpart. One of the main problem is then that the scale quantization
procedure of the Euler-Lagrange equation is not unique. Indeed, we can first quantify the
Lagrangian of the system and then define a quantized Euler-Lagrange equation, or
wa can quantify directly the classical Euler-Lagrange equation. The main point,
proved in the coherence lemma in $\S$.\ref{cohe}, is that these two procedures coincide. \\

The scale quantization procedure being precisely defined, we can specialise it to the
quantum mechanical case. The principal free parameter in the quantization lies in the
order of the regularity of the non differentiable curve, i.e. its H\"older exponent.
Using the Feynman-Hibbs characterization of quantum paths, as well as the Heisenberg inequalities,
we prove in $\S$.\ref{gencond} that the H\"older regularity of a quantum path is $1/2$.
Using this result, we prove in $\S$.\ref{quantnew} that the quantized analogue of the Newtonian
equation of dynamics is a generalized non linear Schr\"odinger equation. This is done
by introducing a wave-function in $\S$.\ref{waveaction}, which is the direct consequence of the complex
nature of the speed, being itself a consequence of the non differentiability of the curve.
Under special assumptions, which can be interpreted, we recover the classical
Schr\"odinger equation.

\section{Scale quantization procedure for classical Lagrangian systems}
\label{scalemap}
\subsection{Classical Lagrangian systems and Euler-Lagrange equation}

In this article, we only discuss classical Lagrangian systems defined as follow :

\begin{defi}
A Lagrangian $L(x,v,t)$ is called classical if it is of the form kinetic energy+potential,
i.e.
\begin{equation}
L(x,v,t)=K(v) + U(x,t) ,
\end{equation}
where $K(v)$ is a quadratic form.
\end{defi}

The basic example for $K(v)$ is the classical kinetic energy of a particule of mass $m$
given by
\begin{equation}
K(v)=\di {1\over 2} m v^2 .
\end{equation}

The dynamics associated to a Lagrangian system is determined by the
{\it Euler-Lagrange equations}.

\begin{defi}
Let $L(x,v,t)$ be a classical Lagrangian system. The Euler-Lagrange equation associated to
$L$ is the following partial differential equation :
$$
\di {d\over dt} \left (
\di {\partial L\over \partial v} (x(t),v(t),t) \right )
=\di {\partial L\over \partial x} (x(t),v(t),t).
\eqno{(EU)}
$$
\end{defi}

We denote by $E$ the mapping associating to $L$ its Euler-Lagrange equation $(EU)$.

\subsection{The scale quantization procedure}

In this section, we define the {\it scale quantization procedure}, which
formalizes Nottale's approach to quantum mechanics. The terminology
suggests that the quantization procedure follows ideas coming from the theory
of the {\it scale relativity} developped by Nottale \cite{no1}.

\subsubsection{The scale quantization map}

We define a map $Q$ which acts on differential operators, variables and functions.\\

The classical variables $x$, $v$, $t$ have quantized analogues which are denoted by
$X=Q(x)$, $V=Q(v)$ and $T=Q(t)$.\\

In the following, we assume that :

\begin{ass}
We have $Q(t)=t$.
\end{ass}

The time variable has then a specific role, being the only variable not affected by
the quantization procedure. \\

We denote also by
\begin{equation}
X(t)=Q(x(t)),\ V(t)=Q(v(t)) ,
\end{equation}
the quantized version of the position trajectory and speed. \\

The main point is that we don't know for the moment the regularity of $X(t)$ or
$V(t)$.\\

The first algebraic properties of $Q$ is the following :

\begin{proper}[Quantization of maps]
\label{pro4}
We consider a map $L:(x,v,t)\mapsto L(x,v,t)$. The quantized map $Q(L)={\cal L}$ is defined by
\begin{equation}
{\cal L} : (X ,V,t)  \mapsto  L(X,V,t) .
\end{equation}
\end{proper}

As a consequence, if $L$ is differentiable with respect to the variable $x$, $v$ or $t$ then
${\cal L}$ is differentiable with respect to $X$, $V$ or $t$.\\

In order to use $Q$ on differential equations, we must precise its behaviour with
respect to differential operators.

\begin{proper}[Operator]
\label{prop2}
We consider a map of the form $f(t)=L(x(t),v(t),t)$, where $x(t)$ and $v(t)$ are differentiable
functions. The differential operator $d/dt$ acts on operator $f$. By the map $Q$, we
define a quantized operator $Q(d/dt)$ such that\\

i) $Q (\di {d f/dt})=Q (d/dt).Q (f)$,\\

where $Q (d/dt)$ is on operator acting on $Q(f)$, depending on the
regularity of $Q(f)$ with respect to $t$ : \\

*) If $Q(f)(t)$ is differentiable with respect to $t$, then $Q(d/dt)=d/dt$.\\

**) If $Q(f)(t)$ is non differentiable with respect to $t$, then $Q(d/dt)=\Box_{\epsilon} /
\Box t$, where $\epsilon (X,h) >\epsilon >0$, $h$ being a constant.
\end{proper}

The constant $h$ must be fixed by physical constraint. In the following, we consider
$h$ has a free parameter.\\

As $f(t)=L(x(t),v(t),t)$, we have $Q(f)(t)=L(X(t),V(t),t)$. Hence, the regularity of
$Q(f)$ with respect to $t$ depends on the regularity of $X(t)$ and $V(t)$ with respect to
$t$. \\

Moreover, as $v=dx/dt$, we have
\begin{equation}
V(t)=\di Q\left (
\di {d\over dt} \right )
[X(t)] .
\end{equation}
Hence, the regularity of $X(t)$ can induce a change in the form of the speed $V(t)$.

\subsubsection{Scale quantization of the Euler-Lagrange equation}

By the quantization procedure we give the quantized version of the Euler-Lagrange
formula (\ref{qel1}).

\begin{lem}
The quantized Euler-Lagrange equation $Q(EU)$ is given by
\begin{equation}
\label{qel3}
\di Q \left ( {d\over dt}\right ) \left [
{\partial {\cal L} \over \partial V} (X(t),V(t),t) \right ] =
\di {\partial {\cal L} \over \partial X } .
\end{equation}
\end{lem}

{\it Proof}. The action of $Q$ on the classical Euler-Lagrange equation (EU) gives
\begin{equation}
\label{qel1}
\di Q \left ( {d\over dt}\right ) \left [ Q \left (
{\partial L \over \partial v} (x(t),v(t),t) \right ) \right ] =
Q \left ( \di {\partial L\over \partial x} \right ) .
\end{equation}

As $L$ is assumed to be differentiable with respect
to the variables $v$ and $x$, we have, using property \ref{pro4} :

\begin{equation}
\label{qel2}
\di Q \left ( {d\over dt}\right ) \left [
{\partial Q(L) \over \partial Q(v)} (Q(x)(t),Q(v)(t),t) \right ] =
\di {\partial Q(L)\over \partial Q(x) } ,
\end{equation}
With our notations, equation (\ref{qel2}) gives equation (\ref{qel3}). This concludes
the proof of the lemma. $\Box$\\

As a consequence, in order to precise the quantization procedure, we only have to
precise the regularity of $Q(x(t))=X(t)$.

\section{Generic trajectories of Quantum mechanics}
\label{gencond}
In order to precise the quantization procedure, we investigate the regularity of
quantum-mechanical path.

\subsection{Feynman and Hibbs genericity condition}

Feynman and Hibbs have already noted in (\cite{fh},p.176-177) that typical path of quantum-mechanical
particle is continuous and nondifferentiable. More precisely, there exist a
quadratic velocity, i.e. if $X(t)$ denotes the particle trajectory, then

$$
\lim_{t\rightarrow t'} \di {(X(t)-X(t'))^2 \over t-t'} \ \mbox{\rm exists}.
\eqno{(FH)}
$$
As a consequence, we have the following result.

\begin{lem}
\label{fhc}
Under Feynman-Hibbs characterization (FH), we have $X(t) \in H^{1/2}$.
\end{lem}

We can deduce the following result on the Hausdorff dimension of typical paths
of quantum mechanics :

\begin{cor}
Under Feyman-Hibbs characterization (FH), the Hausdorff dimension of $X$ is $1/2$.
\end{cor}

{\it Proof}. As $X\in H^{1/2}$, this follows from (theorem 20.6,\cite{tri},p.310). $\Box$\\

In the contrary, the fractal (or Minkowski-Bouligand) dimension $\Delta$ is given by
\begin{equation}
\Delta (X)=2-(1/2) =3/2 ,
\end{equation}
using (\cite{tri},p.154-155).\\

This result has been discussed in great details by Abbott and Wise \cite{aw}.

\subsection{Heisenberg uncertainty principle}

The non-differentiable character of typical paths of quantum mechanics can be seen as
a consequence of {\it Heisenberg uncertainty relations}. We refer to (\cite{no1},p.93-95) and
(\cite{cospi},p.130-131) for details.\\

Let $\Delta x$, $\Delta t$ and $\Delta p$ be the precision of the measurement of the
position $x$, time $t$ and $momentum$ $p$ of a given particle. The Heisenberg
uncertainty relation on momentum and position is
\begin{equation}
\label{hr2}
\Delta p \Delta x \geq h .
\end{equation}
We have the following relations
\begin{equation}
\left .
\begin{array}{lll}
x(t\pm \Delta t) & = & x(t)\pm \Delta x , \\
v_{\epsilon} (t\pm \Delta t) & = & (x(t+\epsilon )-x(t)/\epsilon )
\pm\Delta v ,
\end{array}
\right .
\end{equation}
where $\epsilon \geq \Delta t$ by definition of $\Delta t$, and
\begin{equation}
\Delta p = m \Delta v .
\end{equation}

\begin{rema}
The speed of a particle $v(t)$ is defined as the limit of the mean-speed
$v_{\epsilon} (t) =(x(t+\epsilon )-x(t) /\epsilon )$ when $\epsilon$ goes to zero.
As the time variable $t$ is known only with precision $\Delta t$, we must
have $\epsilon \geq \Delta t$.
\end{rema}

What are the relations between $\Delta v$, $\Delta x$ and $\Delta t$ ? \\

As $\epsilon \geq \Delta t$ by assumption, we have by taking the best
possible value $\epsilon =\Delta t $,
\begin{equation}
\Delta v \sim 2\Delta x /\Delta t .
\end{equation}
As a consequence, the Heisenberg uncertainty relation (\ref{hr2}) gives
\begin{equation}
m \di {(\Delta x)^2 \over \Delta t} \sim h .
\end{equation}
We deduce that
\begin{equation}
\Delta x \sim \di {h\over 2m} \Delta t^{1/2} .
\end{equation}
Hence, we deduce that Heisenberg uncertainty relation (\ref{hr2}) induce the fact that
$X\in H^{1/2}$. \\

Of course, the previous reasoning can be reversed. If typical paths of
quantum mechanics are assumed to be in $H^{1/2}$ then we obtain Heisenberg-like uncertainty
relations.

\begin{lem}
Let $X\in H^{1/2}$, and $0<\Delta t <<1$ be a small parameter. We denote by ${\cal X}$ the
graph of $X$, and by ${\cal X}(t)=(t,X(t))$ a point belonging to ${\cal X} \in \rR^2$. We
denote by $\Delta x (t)=\parallel {\cal X}(t), {\cal X}(t+\Delta t) \parallel$, where
$\parallel .\parallel$ is the classical Euclidean norm on $\rR^2$. We have
\begin{equation}
\Delta x \sim (\Delta t)^{1/2} .
\end{equation}
for all $t\in \rR$.
\end{lem}

{\it Proof}. This follows from a simple computation. We have
\begin{equation}
\parallel {\cal X}(t) ,{\cal X} (t+\Delta t) \parallel ^2 =
(\Delta t)^2 + (X(t+\Delta t) -X(t))^2 .
\end{equation}
As $X\in H^{1/2}$, we obtain
\begin{equation}
\parallel {\cal X}(t) ,{\cal X} (t+\Delta t) \parallel ^2
\leq (\Delta t)^2 +C^2 \Delta t \sim \Delta t ,
\end{equation}
for $\Delta t$ sufficiently small. $\Box$\\

Of course, this result extends to arbitrary $X\in H^{\alpha}$, with $0<\alpha <1$, for
which we obtain
\begin{equation}
\Delta x \sim ( \Delta t )^{\alpha } .
\end{equation}

\section{Scale quantization of Newtonian mechanics}

\subsection{Quantization of speed}

The fact that typical path of quantum mechanics belongs to $H^{1/2}$ has
many implications with respect to the quantization procedure. The first one
being that the quantized version of the speed $v$ is now complex.

\begin{lem}[Quantized speed]
Let $Q(x)=X \in H^{1/2}$, then $Q(v)=V \in \cC$. Precisely, we have
\begin{equation}
V=\di {\Box_{\epsilon} X\over \Box t} .
\end{equation}
\end{lem}

{\it Proof}. This follows from the quantization procedure. We have
\begin{equation}
v=\di {dx\over dt} .
\end{equation}
By $Q$, we obtain
\begin{equation}
V=Q\left (
\di {d\over dt}
\right )
(X).
\end{equation}
As $X\in H^{1/2}$, we have
\begin{equation}
Q(d/dt)=\Box_{\epsilon} /\Box t ,
\end{equation}
which is a complex number by definition. $\Box$

\subsection{Scale relativity and the coherence lemma}
\label{cohe}

\subsubsection{Scale relativity and Scale Euler-Lagrange equations}

The scale-relativity theory developped by Nottale (\cite{no1},\cite{no2})
extends the Einstein relativity principle to {\it scale}.
An heuristic version of the new relativity principle can be written as :\\

``{\it The equations of physics keep the same form under any transformation of scale
(contractions and dilatations)}."\\

The mathematical foundation of such a theory is difficult. The main difficulty being that
space-time is now a non differentiable manifold. We refer to the work of Nottale \cite{no1}
for more details. \\

The scale relativity principle has a direct consequence on the form of the
equation of motion for a particle.\\

We first introduce the following ``Scale" Euler-Lagrange equation :

\begin{defi}
Let ${\cal L}(V,t)$ be a quantized Lagrangian. The Scale Euler-Lagrange equation
associated to ${\cal L}(X,V,t)$ is the equation
$$
\di {\Box_{\epsilon} \over \Box t} \left (
\di {\partial {\cal L} \over \partial V}(X(t),V(t),t)
\right )
=\di {\partial {\cal L} \over \partial X} (X(t),V(t),t).
\eqno{({\cal EU})}
$$
\end{defi}

We denote by ${\cal E}$ the mapping associating equation $({\cal EU})$ to ${\cal L}$.\\

The scale relativity principle is then equivalent to the following statement :

\begin{sta}[Scale relativity]
The equation of motion of a quantum-mechanical particle satisfy the
Scale Euler-Lagrange equation.
\end{sta}

The main point is that the Scale Euler-Lagrange equation doesn't follow from the
quantization procedure, but from a first principle which fixes the form of the equation
of motions.

\subsubsection{Coherence lemma}

At this point, we can state the {\it coherence lemma}, which ensures us that the
quantization procedure is well defined. Indeed, the scale relativity principle gives
an equation of motion particle in quantum mechanics, i.e. $({\cal EU})$. But,
the quantization procedure can be used to obtain an equation of motion for the particle as a quantized version
of the classical Euler-Lagrange equation, i.e. $Q(EU)$. The main point is that this two
constructions are equivalent, i.e. $Q(EU)={\cal EU}$. Precisely, we have :

\begin{lem}[Coherence]
The following diagramm commutes
\begin{equation}
\left .
\begin{array}{ccc}
L(x,v,t) & \stackrel{Q}{\longrightarrow} & {\cal L} (X,V,t) \\
E \downarrow & & \downarrow {\cal E} \\
\di {d\over dt} \left [
\di {\partial L\over \partial v} (x(t),v(t),t) \right ] =
\di {\partial U\over \partial x} & \stackrel{Q}{\longrightarrow} &
\di {\Box_{\epsilon} \over \Box t}
\left [
\di {\partial {\cal L}\over \partial V} (X(t),V(t),t) \right ] =
\di {\partial U\over \partial X} .
\end{array}
\right .
\end{equation}
\end{lem}

{\it Proof}. This follows from a direct computation. $\Box$\\

In term of mapping, the coherence lemma is then equivalent to
\begin{equation}
Q\circ E = {\cal E} \circ Q .
\end{equation}

\subsection{Action functional and wave}
\label{waveaction}
A basic element of Lagrangian mechanics is the action functional $A(x,t)$ which is
related to speed via the equation
\begin{equation}
\label{sa}
v=\di {1\over m} \di {\partial A \over \partial x} .
\end{equation}
The function $A(x,t)$ is differentiable with respect to $x$ and $t$. By the
quantization procedure, we obtain the analogue of the classical action
functional for quantum mechanics :

\begin{lem}[Quantized action functional]
The function ${\cal A} =Q(A)$ is a complex valued function ${\cal A}(X,t)$ which
satisfies
\begin{equation}
\label{qsa}
V=\di {1\over m} \di {\partial {\cal A} \over \partial X} .
\end{equation}
\end{lem}

{\it Proof}. By $Q$ the classical equation (\ref{sa}) gives
\begin{equation}
\label{qsa1}
V=Q(v)=\di Q\left ( {\partial \over \partial x} \right )
Q(A) .
\end{equation}
As the classical action $A(x,t)$ is differentiable with respect to $x$, the
quantized version ${\cal A}(X,t)$ is differentiable with respect to $X$. As
a consequence, we have
\begin{equation}
Q\left ( {\partial \over \partial x} \right ) =\di {\partial \over \partial X} .
\end{equation}
This concludes the proof. $\Box$\\

Moreover, the quantization procedure gives us the following relation :

\begin{lem}[Action]
The quantized action ${\cal A}(x,t)$ satisfies
\begin{equation}
{\cal L}=\di {\partial {\cal A} \over \partial t} .
\end{equation}
\end{lem}

At this point, the main feature of quantum mechanics is to introduce complex speed and
action functional. We can discuss the behaviour of $\cal A$ by introducing a complex
valued function, which is called the {\it wave function}.

\begin{defi}[Wave function]
We called wave function associated to $X$, the complex valued function defined by
\begin{equation}
\psi (X,t)=\di \exp \left (
\di { i {\cal A} (X,t) \over 2m\gamma }
\right ) ,
\end{equation}
where $\gamma \in \rR$ is a constant number.
\end{defi}

The constant $\gamma$ is a normalization constant, which depends on the regularity
property of $X(t)$. \\

The complex nature of $V$ leads naturally to the introduction of the wave function. The
wave formalism is then induced by the non differentiable character of typical paths of
quantum mechanics.

\subsection{The quantized Euler-Lagrange equation}
\label{quantnew}
The quantized Euler-Lagrange equation is given by
\begin{equation}
\label{qele}
m\di {\Box_{\epsilon } V(t)\over  \Box t} =\di {dU \over dx} (x) .
\end{equation}
The complex speed $V$ is related to the wave function, so that equation
(\ref{qele}) can be writen in term of $\psi$.

\begin{thm}
Let $X(t)$ be a continuous non differentiable function in $H^{1/2}$ and $\psi $ its
associated wave function. The quantized Euler-Lagrange equation is of the form
\begin{equation}
\label{gse}
2i\gamma m
\left [
-\di {1\over \psi} \di \left (
{\partial \psi \over \partial X}\right ) ^2
\left (
\di i\gamma +{a_{\epsilon} (t)\over 2}
\right )
+\di {\partial \psi \over \partial t}
+\di {a_{\epsilon} (t)\over 2} \di {\partial^2 \psi \over \partial X^2} \right ]
=(U(x) +\alpha (x) )\psi +o(\epsilon^{1/2} ),
\end{equation}
where $\alpha (x)$ is an arbitrary continuous function, and
\begin{equation}
\label{aeps}
a_{\epsilon } (t)=\di {1\over 2} \left [
(\Delta_{\epsilon}^+ X(t))^2 -
(\Delta_{\epsilon}^- X(t))^2
\right ] -i \di {1\over 2}
\left [
(\Delta_{\epsilon}^+ X(t))^2 +
(\Delta_{\epsilon}^- X(t))^2
\right ]
.
\end{equation}
\end{thm}

The equation (\ref{gse}) is called {\it generalized Schrodinger equation} by Nottale
\cite{no1}.\\

{\it Proof}. We have
\begin{equation}
V=-i2\gamma \di {\partial \ln (\psi )\over \partial X} .
\end{equation}
The Euler-Lagrange equation is now given by
\begin{equation}
2i\gamma m \di {\Box_{\epsilon} \over \Box t} \di \left (
\di {\partial \ln (\psi ) \over \partial X} \right ) =\di {dU\over dX} .
\end{equation}
We denote
\begin{equation}
f(X,t)=\di {\partial \ln (\psi (X,t))\over \partial X} (X,t) .
\end{equation}
We apply the main lemma of part I to compute $\Box_{\epsilon} f(X(t),t) /\Box t$ for
complex valued functions. We have
\begin{equation}
\label{cal1}
\left .
\begin{array}{lll}
\di {\Box_{\epsilon} \over \Box t} \di \left (
\di {\partial \ln (\psi ) \over \partial X} (X(t),t) \right )
& = & \di {\Box_{\epsilon} X\over \Box t} \di {\partial \over \partial X} \left (
\di {\partial \ln (\psi (X,t))\over \partial X}
\right ) (X(t),t) \\
 & & +
\di {\partial \over \partial t}
\left (
\di {\partial \ln (\psi (X,t))\over \partial X}
\right ) (X(t),t) \\
 & & +\di {1\over 2} a_{\epsilon} (t) \di {\partial^2 \over \partial X^2}
\left (
\di {\partial \ln (\psi (X,t))\over \partial X}
\right ) (X(t),t)
+o(\epsilon^{1/2} ) .
\end{array}
\right .
\end{equation}
Elementary calculus gives
\begin{equation}
\di {\partial \ln (\psi (X,t))\over \partial X} =\di {1\over \psi} \di {\partial \psi\over
\partial X},\ \ \mbox{\rm and}\ \ \
\di {\partial \over \partial X}
\left (
{1\over \psi} \di {\partial \psi\over
\partial X}
\right ) =
{1\over \psi} \di {\partial^2 \psi\over
\partial^2 X}
-
{1\over \psi^2} \di \left ( {\partial \psi\over
\partial X}\right ) ^2 .
\end{equation}
Moreover, by definition of the wave function $\psi$ and $V$, we have
\begin{equation}
V=-i2\gamma \di {\partial \ln \psi (X,t) \over \partial X} .
\end{equation}
Hence, we obtain
\begin{equation}
\di {\Box_{\epsilon} X\over \Box t} =V=
-i2\gamma \di {\partial \ln \psi (X,t) \over \partial X} ,
\end{equation}
and
\begin{equation}
\left .
\begin{array}{lll}
\di {\Box_{\epsilon} X\over \Box t} \di {\partial \over \partial X} \left (
\di {\partial \ln (\psi (X,t))\over \partial X}
\right ) (X(t),t) & = &
-i2\gamma \di {\partial \ln (\psi )\over \partial X}
\di {\partial \over \partial X}
\left (
\di {\partial \ln (\psi )\over \partial X} \right )
(X(t),t) , \\
 & = & -i\gamma \di {\partial \over \partial X}
 \left [
 \di \left (
 \di {\partial \ln (\psi )\over \partial X}
 \right )
 ^2
 \right ]
(X(t),t) , \\
 & = & -i\gamma \di {\partial \over \partial X}
 \left [
 \di {1\over \psi^2}
 \di \left (
 \di {\partial \psi \over \partial X}
\right ) ^2
\right ]
(X(t),t) .
\end{array}
\right .
\end{equation}
We then have
\begin{equation}
\left .
\begin{array}{l}
\di {\Box_{\epsilon} \over \Box t}
\left (
\di {\partial \ln (\psi (X,t))\over \partial X} (X(t),t)
\right )  \\
= \di {\partial \over \partial X}
\left [
-i\gamma \di {1\over \psi^2}
\left (
\di {\partial \psi \over \partial X}
\right ) ^2
+\di {\partial \ln (\psi ) \over \partial t}
  +\di {1\over 2} a_{\epsilon} (t)
\left [
\di {1\over \psi}
\di {\partial^2 \psi \over \partial X^2} -\di {1\over \psi^2}
\left (
\di {\partial \psi \over \partial X}
\right ) ^2
\right ]
\right ]
+o(\epsilon^{1/2} ) ,\\
= \di {\partial \over \partial X} \left [
 -{1\over \psi^2}
\left ( \di {\partial \psi\over \partial X} \right ) ^2
\left (
i\gamma +\di {a_{\epsilon} (t) \over 2} \right )
+\di {1\over \psi} \di {\partial \psi\over\partial t}
+\di {a_{\epsilon} (t) \over 2} \di {1\over \psi} \di
{\partial^2 \psi \over \partial X^2}
\right ] +o(\epsilon^{1/2} ).
\end{array}
\right .
\end{equation}

As a consequence, equation (\ref{cal1}) is equivalent to
\begin{equation}
\di {\partial\over \partial X}
\left [
i2\gamma m
\left [
-{\di 1\over \psi^2} \left (
\di {\partial \psi \over \partial X}
\right )
^2
\left (
i\gamma  +\di {a_{\epsilon} (t) \over 2}
\right )
+\di {1\over \psi} \di {\partial \psi\over\partial t} \right ]
+\di {a_{\epsilon} (t) \over 2} \di {1\over \psi} \di {\partial^2 \psi \over
\partial X^2}
\right ] =\di {\partial U\over \partial X}.
\end{equation}
By integrating with respect to $X$, we obtain
\begin{equation}
i2\gamma m
\left [
-{\di 1\over \psi^2} \left (
\di {\partial \psi \over \partial X}
\right )
^2
\left (
i\gamma  +\di {a_{\epsilon} (t) \over 2}
\right )
+\di {1\over \psi} \di {\partial \psi\over\partial t} \right ]
+\di {a_{\epsilon} (t) \over 2} \di {1\over \psi} \di {\partial^2 \psi \over
\partial X^2}
=U(X)+\alpha (X)+o(\epsilon^{1/2} ),
\end{equation}
where $\alpha (X)$ is an arbitrary function. This concludes the proof. $\Box$

\subsection{Nonlinear Schr\"odinger equations}

Many authors have suggested that the quantum mechanics based on the linear
Schr\"odinger equation is only an approximation of some non linear theory with a
nonlinear Schr\"odinger equation. Nonlinear wave mechanics was initiated by
De Broglie \cite{debr} in order to have a better understanding of the relation
between wave and particle (see \cite{debr2},p.227-231).\\

Many generalizations of the Schr\"odinger equation exist. We can mention for example
the {\it Staruszkievicz} \cite{puz} or {\it Bialynicki-Birula and Mycielsky} \cite{bbm} modification.
All these generalizations are not a consequence of a given principle of physic. For
example, the Bialynicki-Birula and Mycielsky modification can be derived, using the
hydrodynamical formalism proposed by Madelung \cite{mad}, Bohm-Vigier \cite{bv} and others,
by adding a {\it pressure term} to the Euler hydrodynamical equation (see \cite{pard},$\S$.2).
Of course, one can justify a posteriori a given modification by the fact that it solves some
relevant problems like the {\it collaps} of the wave function or the Schr\"odinger cat
paradox (see \cite{pard}).\\

Recently Castro, Mahecha and Rodriguez \cite{cmr} have proposed a new nonlinear
Schr\"odinger equation based on a generalization of Nottale reasoning \cite{no1}. In the
one-dimensional case, they obtain
\begin{equation}
\label{nngse}
i\alpha \di {\partial \psi\over \partial t}
=-\di {\alpha \re (\alpha ) \over 2m}
\di {\partial^2 \psi \over \partial x^2}
+U(x) \psi -i\di {\alpha \im (\alpha ) \over 2m}
\left (
\di {\partial \psi \over \partial x}
\right )^2 \psi ,
\end{equation}
where $\alpha \in \cC$.\\

If $\im (\alpha )=0$ and $\alpha =\bar{h}$, then one recover the classical linear
Schr\"odinger equation.\\

This generalization is done under two assumptions :\\

i) They introduce a complex ``Planck constant" (i.e. $\alpha$) by allowing that
the normalization constant $\gamma$ be complex;\\

ii) They consider a {\it complex diffusion coefficient}, which has nos counterpart in
our case.\\

The problem is that there is no geometric interpretation of a
complex diffusion coefficient, nor of a complex Planck constant. As a consequence,
even if Nottale's reasoning is keep in order to derive the new non linear
Schr\"odinger equation (\ref{nngse}), this is an ad-hoc mathematical generalization.\\

In the contrary, we have obtain a non linear Schr\"odinger equation (\ref{gse})
directly from the Scale relativity principle. All our constants have a clear
geometrical meaning. Moreover, the nonlinear term is not an ad-hoc term but fixed
by the theory.

\section{Toward the Schr\"odinger equation}

The generalized Schr\"odinger equation can be simplified in some case, by assuming that
quantum-mechanical paths satisfy special regularity properties.

\begin{thm}
Let $X(t)$ be a continuous non differentiable function belonging to $H^{1/2}$ such that
\begin{equation}
\label{regcond}
a_{\epsilon} (t)=-i2\gamma .
\end{equation}
Then, the quantized Euler-Lagrange equation takes the form
\begin{equation}
\label{mschro}
\gamma^2 \di {\partial^2 \psi \over \partial X^2} +i\gamma \di {\partial \psi
\over \partial t} =
(U(X,t)+\alpha (X) )\di {\psi \over 2m} +o(\epsilon^{1/2}).
\end{equation}
We can always choose a solution of equation (\ref{mschro}) such that
\begin{equation}
\label{phase}
\alpha (X)=0.
\end{equation}
In this case, if
$$
\gamma =\di {\bar{h} \over 2m} ,
\eqno{(SC)}
$$
we obtain the classical Schr\"odinger equation
\begin{equation}
i\bar{h} \di {\partial \psi\over \partial t} +\di {\bar{h}^2 \over 2m} \di {\partial ^2
\psi \over \partial X^2} =U\psi .
\end{equation}
\end{thm}

{\it Proof}. The only non trivial, but classical, part concerns the possibility to
choose a phase of the wave function such that $\alpha (X)=0$. The computations
follow closely those of (\cite{bc1},corollary 1).\\

Let $\psi$ be a solution of (\ref{mschro}). The basic idea is that we can always modify the
phase of the wave function in order to obtain a solution of (\ref{mschro}) such that
$\alpha (X)=0$. So, let us consider the modified wave function
\begin{equation}
\tilde{\psi} (X,t)=\di \exp \left (
i\di {{\cal A} (X,t) \over 2m\gamma} +\theta (X) \right ) =\psi (X,t) \Theta (X).
\end{equation}
We have
\begin{equation}
\left .
\begin{array}{lll}
\di {\partial \tilde{\psi} \over \partial X} & = & \di {\partial \psi \over \partial X} \Theta
+\psi \Theta ' ,\\
\di {\partial^2 \tilde{\psi} \over \partial X^2} & = & \di {\partial^2 \psi \over \partial X^2}
\Theta +2 {\partial \psi \over \partial X} \Theta ' +\psi \Theta " ,\\
\di {\partial \tilde{\psi} \over \partial t} & = &
\di {\partial \psi \over \partial t} \Theta ,
\end{array}
\right .
\end{equation}
where $\Theta '(X)$ and $\Theta"(X)$ are the first and second derivative of $\Theta$.\\

By replacing in (\ref{mschro}), and assuming that $\alpha (X)=0$, we obtain up to
$o(\epsilon^{1/2})$ terms,
\begin{equation}
\Theta \left (
i2\gamma m \di {\partial \psi \over \partial t} +2\gamma^2 m \di {\partial^2 \over \partial
X^2} -U\psi \right )
+4\gamma^2 m \di {\partial \psi \over \partial X} \Theta' +2\gamma^2 m\psi \Theta " =0 .
\end{equation}

As $\psi (X,t)$ is a solution of (\ref{mschro}) with a given value $\alpha (X)$, we deduce
that $\Theta (X)$ must satisfy the following ordinary differential equation
\begin{equation}
a(X,t) \Theta  +b(X,t) \Theta ' +c(X,t) \Theta " =0,
\end{equation}
where
\begin{equation}
a(X,t)=\alpha (X) \psi , \ b(X,t)=4\gamma^2 m  \di {\partial \psi \over \partial X} , \
c(X,t)=2\gamma^2 m \psi .
\end{equation}
This is a second order differential equation with non constant coefficients. By
general theorems on linear differential equations, there always exists a solution. As a consequence, we can always choose
$\Theta (X)$ such that $\phi(X,t) \Theta (X)$ satisfies (\ref{mschro}) with $\alpha (X)=0$.
This concludes the proof. $\Box$

\subsection{Difference equations and the Schr\"odinger condition}

The Schr\"odinger condition (SC) can be precised. We have the following lemma :

\begin{lem}
The Schr\"odinger condition (SC) is equivalent to the difference equation
\begin{equation}
\Delta_{\epsilon}^+ X(t)=\pm\sqrt{\bar{h} /m} ,
\end{equation}
and the following relations
$$
\Delta_{\epsilon}^+ X(t)=\Delta_{\epsilon}^- X(t) . \eqno{(P)}
$$
\end{lem}

{\it Proof}. By definition of $a_{\epsilon} (t)$ in (\ref{aeps}), the Schrodinger condition is
equivalent to the following system
\begin{equation}
\left .
\begin{array}{lll}
(\Delta_{\epsilon}^+ X(t))^2 -(\Delta_{\epsilon}^- X(t))^2 & = & 0 , \\
(\Delta_{\epsilon}^+ X(t))^2 +(\Delta_{\epsilon}^- X(t))^2 & = & 2\bar{h}/m .
\end{array}
\right .
\end{equation}
We deduce that $(\Delta_{\epsilon}^+ X(t))^2 =(\Delta_{\epsilon}^- X(t))^2$ and
$(\Delta_{\epsilon}^+ X(t))^2 =\bar{h}/m$. Hence, we have
$\Delta_{\epsilon}^+ X(t)=\pm \sqrt{\bar{h} /m}$, a constant independant of $t$. We
obtain
\begin{equation}
\Delta_{\epsilon}^- X(t)=\Delta_{\epsilon}^+ (t-\epsilon )=
\pm \sqrt{2\bar{h} /m} =\Delta_{\epsilon}^+ X(t) ,
\end{equation}
which concludes the proof. $\Box$\\

As a consequence, we are lead to the study of difference equations of the form
\begin{equation}
\label{modequ}
\Delta_{\epsilon}^{\sigma} X(t)=a, \ \sigma =\pm ,
\end{equation}
where $a\in\rR$ is a {\it constant}. This kind of difference equations\footnote{
Most of the following can be generalized to general difference equations of the
form $\Delta_{\epsilon}^{\sigma} X(t)=\Phi (t)$, where $\sigma =\pm$ and $\Phi (t)$ is a given
function (see \cite{mil},Chap.8).} always have
solutions of the form
\begin{equation}
X(t)=X^* (t)+P_{\epsilon} (t) ,
\end{equation}
where $X^* (t)$ is a particular solution, and $P_{\epsilon} (t)$ is an {\it arbitrary
periodic} function of $t$ of period $\epsilon$.\\

Particular solutions always exists. Indeed, if one consider an arbitrary given function
$X^* (t)$ defined on $0\leq t<\epsilon$, then the difference equation define $X^* (t)$
at every point exterior to this interval. Of course, such kind of solutions are in
general not {\it analytic}. We can define a special particular solution called
{\it principal solution} following (\cite{mil},p.200).\\

For (\ref{modequ}) the principal solution is given by (see \cite{mil},p.204)
\begin{equation}
X^*_c (t)=a(t-c-\di {\epsilon \over 2} ) ,
\end{equation}
where $c$ is an arbitrary constant.

\begin{rema}
For $c=0$, $\omega=1$, $a=1$, we obtain the Bernouilli's polynomial $B_1 (t)$.
\end{rema}

>From this section, we deduce the following {\it structure lemma} :

\begin{lem}[structure]
The Schr\"odinger condition is satisfied by continuous functions
of the form
\begin{equation}
X_c (t)=\pm \sqrt{\bar{h}/m} \,  \left (  t-c- \di {\epsilon \over 2} \right ) +P_{\epsilon} (t),
\end{equation}
where $c$ is an arbitrary real constant, and $P_{\epsilon} (t)$ is an
arbitrary periodic function belonging to $H^{1/2}$.
\end{lem}

We call this set of functions, the {\it principal Schr\"odinger set}. It gives
the structure of typical paths of quantum mechanics in the free case.

\begin{rema}
In \cite{bc1}, the Schr\"odinger equation is obtained from a generalized equation
under an assumption which is
equivalent to a special fractional differential equation. In \cite{bc4}, we have prove
that in the framework of fractional calculus, this fractional equation has no solutions.
\end{rema}

\end{document}